\documentclass[11pt,reqno]{amsproc}
\usepackage{cite}
\usepackage{fullpage}
\usepackage{graphicx}
\usepackage{subfigure}
\usepackage{psfrag}
\usepackage{color}
\usepackage{amsmath}
\usepackage{amssymb}
\usepackage{enumerate}

\usepackage{algorithm}
\usepackage{algorithmic}

\title{A mixed formulation for a modification to Darcy equation based on Picard linearization and numerical solutions to large-scale realistic problems}

\author{K.~B.~Nakshatrala} \author{D.~Z.~Turner}

\address{Correspondence to: Dr.~Kalyana Babu Nakshatrala, Department 
of Civil and Environmental Engineering, Engineering Building \#1, Room 
\#135,  University of Houston, Houston, Texas - 77204. 
TEL: +1-713-743-4418} \email{knakshatrala@uh.edu}

\address{Dr.~Daniel Z.~Turner, Department of Civil Engineering, 
University of Stollenbosch, Stollenbosch, South Africa. 
TEL: +27-21-808-4434} \email{dzturner@sun.ac.za}

\date{\today}

\begin{document}

\begin{abstract}
  In this paper we consider a modification to Darcy equation by taking 
  into account the dependence of viscosity on the pressure. We present 
  a stabilized mixed formulation for the resulting governing equations. 
  Equal-order interpolation for the velocity and pressure is considered, 
  and shown to be stable (which is not the case under the classical 
  mixed formulation). The proposed mixed formulation is tested 
  using a wide variety of numerical examples. The proposed 
  formulation is also implemented in a parallel setting, and the 
  performance of the formulation for large-scale problems is 
  illustrated using a representative problem. Two practical and 
  technologically important problems, one each on enhanced 
  oil recovery and geological carbon-dioxide sequestration, 
  are solved using the proposed formulation. The numerical 
  examples show that the predictions based on Darcy model 
  are qualitatively and quantitatively different from that of the 
  predictions based on the modified Darcy model, which takes 
  into account the dependence of the viscosity on the pressure. 
  In particular, the numerical example on the geological carbon-dioxide 
  sequestration shows that Darcy model over-predicts the leakage into 
  an abandoned well when compared to that of the modified Darcy model. 
  On the other hand, the modified Darcy model predicts higher pressures 
  and higher pressure gradients near the injection well. These predictions 
  have dire consequences in predicting damage and fracture zones, and 
  designing the seal, whose integrity is crucial to the safety of a geological 
  carbon-dioxide sequestration geosystem. 
\end{abstract}

\maketitle

\section{INTRODUCTION}
\label{Sec:MDarcy_Introduction}
Darcy equation has been successfully employed to model flow in porous 
media in wide areas of applications ranging from groundwater hydrology, 
petroleum engineering to food engineering. Henry Darcy has proposed a 
simple equation to model the flow of an incompressible fluid in rigid 
porous media, which is popularly referred to as Darcy equation and in 
some instances as ``Darcy law.'' Darcy has developed this equation 
empirically based on his experiments on the flow of water through 
sand beds \cite{Darcy_1856}. (Also see the English translation of 
Darcy's work by Patricia Bobeck \cite{Bobeck_Darcy}.) 

It is important to note that Darcy equation is just an approximation of 
the balance of linear momentum for the fluid flowing through a rigid porous 
solid (see the discussion in \cite[Section 2.5]{Bowen_Porous_elasticity}). 
The equation is valid under a number of assumptions, see references 
\cite{Bowen_Porous_elasticity,Rajagopal_M3AS_2007_v17_p215} and also 
\cite[Introduction]{Nakshatrala_Rajagopal_IJNMF_arxiv_2010}. Darcy equation 
merely predicts the flux, and this prediction is \emph{not} accurate at 
high pressures. Moreover, in references 
\cite{Munaf_Wineman_Rajagopal_Lee_M3AS_1993_v3_p231,Nakshatrala_Rajagopal_IJNMF_arxiv_2010} 
it has been advocated that Darcy equation is quite good for low flow 
rates for a fluid like water but is not the case for high flow rates 
and for dense fluids like mineral oil.  

Several extensions to Darcy equation have been developed by 
various researchers. Two of the early popular ones are by 
Forchheimer \cite{Forchheimer_1901_v45_p1782} and Brinkman 
\cite{Brinkman_ASR_1947_vA1_p81}. Bowen has outlined various 
models for flow through porous media allowing multiple 
fluid components and deformation of the porous solid 
\cite{Bowen_Porous_elasticity}. Recent extensions of 
Darcy equation can be found in references 
\cite{2006_Medeiros_Gurgel_Marcondes_JPM_v9_p235,
2011_Farhany_Turan_IJHMT_v54_p2868,
2012_Mahapatra_Pal_Mondal_NAMC_v17_p223}. 
A hierarchy of models for flow through porous media has been 
presented in \cite{Rajagopal_M3AS_2007_v17_p215} of which the 
simplest model is Darcy equation. One important point to note, 
which is central to this paper, is that Darcy equation assumes 
that the coefficient of viscosity (and hence the drag coefficient) 
to be independent of the pressure in the fluid. 

\subsection{Pressure dependent viscosity}
Certain fluid models do involve flow equations with 
non-constant viscosity (e.g., a function of pressure 
or velocity), which are more realistic in modeling 
enhanced oil recovery and geological carbon-dioxide 
sequestration. In both enhanced oil recovery and 
geological carbon-dioxide sequestration, pressure can 
vary from $0.1$ MPa to $100$ MPa. There is irrefutable 
experimental evidence that viscosity of mineral oils 
is not constant and changes drastically with respect 
to pressure (for example; Bridgman \cite{Bridgman}, 
Andrade \cite{Andrade_PRSLA_1952_v215_p36}). In fact, 
Barus \cite{Barus_AJS_1893_v45_p87} suggested the 
following exponential relationship between viscosity 
and pressure: 
\begin{align}
  \label{Eqn:StabMDarcy_Barus_formula}
  \mu(p) = \mu_0 \exp[\beta p]
\end{align}
where $\beta$ has units $\mathrm{Pa}^{-1}$. 

In geological carbon-dioxide sequestration, supercritical 
carbon-dioxide is pumped into deep saline aquifers or 
abandoned oil wells. In references 
\cite{Vesovic_Wakeham_Olchowy_Sengers_Watson_Millat_JPhysChem_1990_v19_p763,
Abramson_PRE_2009_v80_021201}, transport properties (including viscosity) of 
supercritical carbon-dioxide are presented for various pressures ranging from 
$0.1$ MPa to $8$ GPa, and for various temperatures. These experimental studies 
clearly indicate that the viscosity varies exponentially with respect to pressure 
even for supercritical carbon-dioxide for various temperatures and for a wide 
range of pressures (varying from $0.1$ MPa to $8$ GPa).

It should be emphasized that petroleum reservoir simulations 
and geological carbon-dioxide sequestration involve heat 
transfer, mass transfer and phase changes in addition to 
flow aspects. Specifically, one has take into account 
the attendant geochemical reactions (e.g., see Lichtner 
\cite{Lichtner_RMG_1996_v34_p1}). Several mathematical 
models, numerical formulations, and numerical simulations 
have been reported in the literature that have taken 
into account more than the flow aspects. For example, 
see references \cite{1985_Lichtner_GeCA_1985_v49_p779,
Lu_Litchner_JPCS_2007_v78,2010_Marcondes_Sepehrnoori_JPSE_v73_p99}.
The mathematical model in this paper does not address such issues 
but just focuses on the need to take into account the dependence 
of viscosity on the pressure in enhanced oil recovery and geological 
carbon-dioxide sequestration simulations. 
Our study is primarily motivated by the fact that the flow 
characteristics of a fluid whose viscosity depends on pressure 
can be significantly different from that of the flow characteristics 
of a fluid with constant viscosity. Some representative prior works 
on modeling of fluids by taking into account the dependence of viscosity 
on pressure are \cite{Malek_Necas_Rajagopal_ARMA_2002_v165_p243,
Subramaniam_Rajagopal_CMA_2007_v53_p260,
Bulicek_Malek_Rajagopal_IUMJ_2007_v56_p51,
Kannan_Rajagopal_AMC_2008_v199_p748}. 
\emph{In this paper we consider a modification to Darcy equation 
  that takes into account the dependence of viscosity on pressure. 
  The resulting equations will be in mixed form and nonlinear. 
  The unknowns are the velocity and pressure. We present a new 
  stabilized mixed formulation based on the variational multiscale 
  formalism and fixed-point linearization.}

\subsection{Stabilized mixed finite element formulations}
It is well-known that care should be taken to avoid numerical instabilities 
when dealing with mixed formulations. As discussed by Franca and Hughes 
\cite{Franca_Hughes_CMAME_1988_v69_p89}, a stable mixed formulation either 
meets or circumvents the Ladyzhenskaya-Babu\v ska-Brezzi (LBB) \emph{inp-sup} 
stability condition \cite{Babuska_NumerMath_1971_v16_p322,Brezzi_Fortin}. 
Stabilized methods typically fall in the later category, and the stability 
is achieved through addition of stabilization terms.
But most of the approaches that aim to satisfy LBB condition achieve 
the stability by restricting the interpolation functions for the 
independent variables. For example, P1P0 approach in which the 
velocity unknowns are placed at nodes, and pressure unknowns 
are for each element/cell. Two other notable works to satisfy 
the LBB condition are to use either the Raviart-Thomas spaces 
\cite{Raviart_Thomas_MAFEM_1977_p292} or Brezzi-Douglas-Marini 
spaces \cite{Brezzi_Douglas_Marini_NumerMath_1985_v47_p217}. In 
our paper, we do not take such an approach of placing restrictions 
on the interpolation functions for the independent variables. 
Instead, we circumvent the LBB condition by adding stabilization 
terms in a consistent manner, and allow any combination of 
interpolation functions for the independent variables (in this 
case, the velocity vector and the pressure).  


It should be emphasized that the LBB stability condition is 
applicable even to mixed formulations based on the finite 
volume method. As mentioned earlier, a mixed method has to 
either satisfy the LBB condition in the case of saddle-type 
formulation, or circumvent the LBB condition by adding 
stabilization terms to avoid saddle-type problem. In the 
context of a mixed formulation for Darcy-type equation, 
the LBB stability condition places a restriction on the 
choice of the function spaces for the vector-field variable 
(i.e., the velocity) and the scalar-field variable (i.e., 
the pressure). It is interesting to note that staggered 
grids, which is quite popular in the finite volume literature, 
is a way to satisfy the LBB condition.

Some popular approaches for developing stabilized finite element 
formulations are least-squares \cite{Bochev_Gunzburger,Jiang_LSFEM}, 
Galerkin/least-squares (GLS) 
\cite{Hughes_Franca_Hulbert_CMAME_1989_v73_p173,
Baiocchi_Brezzi_France_CMAME_1993_v105_p125}, finite increment calculus (FIC) 
\cite{Onate_CMAME_2000_v182_p355,Onate_IJNME_2004_v60_p255}, streamline upwind 
Petrov-Galerkin \cite{Brooks_Hughes_CMAME_1982_v32_p199}. Another popular approach 
for developing stabilized formulations is the variational multiscale framework 
\cite{Hughes_CMAME_1995_v127_p387}, which has been employed in this paper. The 
variational multiscale framework has been successfully employed in many studies 
to develop stabilized formulations for a wide variety of problems: Darcy equation 
\cite{Masud_Hughes_CMAME_2002_v191_p4341, Hughes_Masud_Wan_CMAME_2006_v195_p3347,
Nakshatrala_Turner_Hjelmstad_Masud_CMAME_2006_v195_p4036}, Stokes' equation 
\cite{Masud_IJNMF_2007_v54_p665,Turner_Nakshatrala_Hjelmstad_IJNMF_2009_v60_p1291}, 
linearized elasticity \cite{Masud_Xia_JAM_2005_v72_p711,
Nakshatrala_Masud_Hjelmstad_CompMech_2008_v41_p547}, incompressible Navier-Stokes 
\cite{Masud_Khurram_CMAME_2006_v195_p1750,Turner_Nakshatrala_Hjelmstad_IJNMF_2010_v62_p119}, 
Fokker-Planck  \cite{Masud_Bergman_CMAME_2005_v194_p1513}. 

A huge volume of literature is available on mixed methods and 
stabilized formulations, and a thorough discussion on these 
topics is beyond the scope of this paper. Some representative 
papers on stabilized formulations are 
\cite{Babuska_Oden_Lee_CMAME_1977_v11_p175,
Brooks_Hughes_CMAME_1982_v32_p199,Oden_Jacquotte_CMAME_1984_v43_p231,
Hughes_Franca_Balestra_CMAME_1986_v59_p85,Franca_Hughes_CMAME_1988_v69_p89,
Doughlas_Wang_MathComput_1989_v52_p495,Hughes_Franca_Hulbert_CMAME_1989_v73_p173,
Franca_Frey_Hughes_CMAME_1992_v95_p253,Franca_Russo_CMAME_1997_v142_p361,
Masud_Hughes_CMAME_2002_v191_p4341,Burman_Hansbo_NMPDE_2005_v21_p986,
Nakshatrala_Turner_Hjelmstad_Masud_CMAME_2006_v195_p4036,Correa_Loula_CMAME_2008_v197_p1525}. 
Some texts concerning stabilized methods are \cite{Quarternio_Valli_PDE,Roos_Stynes_Tobiska,
Chen_Huan_Ma,Donea_Huerta}.
Some representative works on mixed methods in the context of flow through porous media are 
\cite{Ewing_Russell_Wheeler_CMAME_1984_v47_p73,Chavent_Cohen_Jaffre_CMAME_1984_v47_p93,
Ewing_Heinemann_CMAME_1984_v47_p161,Durlofsky_WRR_1994_v30_p965,
Arbogast_Wheeler_Yotov_SIAMJNA_1997_v34_p828,Bergamaschi_Mantica_Manzini_SIAMJSC_1998_v20_p970,
Masud_Hughes_CMAME_2002_v191_p4341,Brezzi_Hughes_Marini_Masud_SIAMJSC_2005_v22_p119,
Hughes_Masud_Wan_CMAME_2006_v195_p3347,Nakshatrala_Turner_Hjelmstad_Masud_CMAME_2006_v195_p4036,
Burman_Hansbo_JCAM_2007_v198_p35}. It should be noted that none of the aforementioned numerical works 
considered the fluid model considered in this paper. 

\subsection{Main contributions of this paper} 
To the best of our knowledge no other prior studies have systematically studied 
(in a numerical setting) the modified Darcy equation except for Reference 
\cite{Nakshatrala_Rajagopal_IJNMF_arxiv_2010}. But this paper presents an 
alternate and simpler formulation than the one presented in Reference 
\cite{Nakshatrala_Rajagopal_IJNMF_arxiv_2010}. The simplicity arises due to 
the fact that the proposed formulation is obtained by first linearizing the 
governing equations and then applying the variational multiscale formalism 
on the resulting equations. This approach facilitates to extend the proposed 
formulation to other complicated models (e.g., modified Brinkman by taking 
into account the pressure dependent viscosity, Brinkman-Forchheimer) in a 
straightforward manner, which may not be the case with the stabilized mixed 
formulation proposed in Reference \cite{Nakshatrala_Rajagopal_IJNMF_arxiv_2010}). 
The main contributions of this paper can be summarized as follows:
\begin{enumerate}[(a)]
\item We considered a realistic modification to Darcy 
  equation by taking into account the dependence of 
  viscosity on the pressure. This implies that the 
  drag coefficient will depend on the pressure. We 
  presented a new stabilized mixed formulation 
  based on the variational multiscale formalism 
  for the resulting nonlinear equations. 
\item We have shown numerically that \emph{equal-order} 
  interpolation for the velocity and pressure (which is 
  computationally the most convenient) is stable under 
  the proposed formulation. A wide variety of test 
  problems are performed to illustrate the performance 
  of the proposed formulation. 
\item We have implemented the proposed in a parallel setting. 
\item We have solved two practical and technologically important 
  large-scale problems with relevance to enhanced oil recovery 
  and geological carbon-dioxide sequestration. The numerical 
  results have clearly indicated the importance of considering 
  the role of dependence of viscosity on the pressure in these 
  two application areas. 
\end{enumerate}

\subsection{An outline of the paper} 
The remainder of the paper is organized as follows. In 
Section \ref{Sec:S2_StabMDarcy_Governing_equations}, we 
present the modified Darcy equation, and present the 
proposed stabilized mixed finite element formulation, 
which is based on the variational multiscale formalism. 
We also present a numerical solution procedure for 
solving the resulting nonlinear equations. In Section 
\ref{Sec:S4_StabMDarcy_Numerical_Results}, we illustrate 
the performance of the proposed formulation on a wide 
variety of benchmark problems, which are commonly used 
in the literature for testing mixed formulations. In the 
same section we also present the numerical results for 
two large-scale practical problems, which are solved by 
implementing the proposed formulation in a parallel 
setting. Conclusions are drawn in Section 
\ref{Sec:S5_StabMDarcy_Conclusions}. 

%
\section{GOVERNING EQUATIONS AND STABILIZED FORMULATION}
\label{Sec:S2_StabMDarcy_Governing_equations}
Let $\Omega \subset \mathbb{R}^{nd}$ be an open and bounded 
set, where ``$nd$'' denotes the number of spatial dimensions. 
Let $\partial \Omega := \bar{\Omega} - \Omega$ be the boundary 
(where $\bar{\Omega}$ is the set closure of $\Omega$), which is 
assumed to be piecewise smooth. A spatial point in $\bar{\Omega}$ 
is denoted by $\boldsymbol{x}$. The spatial gradient and divergence 
operators are, respectively, denoted as ``$\mathrm{grad}[\cdot]$" 
and ``$\mathrm{div}[\cdot]$". 
Let $\boldsymbol{v} : \Omega \rightarrow \mathbb{R}^{nd}$ 
denote the velocity field, and $p : \Omega \rightarrow 
\mathbb{R}$ denote the pressure field. The boundary is 
divided into two parts, denoted by $\Gamma^{v}$ and 
$\Gamma^{p}$, such that $\Gamma^{v} \cap \Gamma^{p} = 
\emptyset$ and $\Gamma^{v} \cup \Gamma^{p}=\partial 
\Omega$. $\Gamma^{v}$ is the part of the boundary on 
which normal component of the velocity is prescribed, 
and $\Gamma^{p}$ is part of the boundary on which 
pressure is prescribed. 

The modified Darcy equations can be written as  
\begin{subequations}
\label{Eqn:MDarcy_modified_Darcy}
\begin{align}
    \label{Eqn:MDarcy_Equilibrium}
    &\alpha(p) \boldsymbol{v} + \mathrm{grad} [p] = 
    \rho(\boldsymbol{x}) \boldsymbol{b}(\boldsymbol{x}) 
    \quad \mbox{in} \; \Omega \\
    \label{Eqn:MDarcy_Continuity}
    &\mathrm{div}[\boldsymbol{v}] = 0 \quad 
    \mbox{in} \; \Omega \\
    \label{Eqn:MDarcy_Velocity_BC}
    &\boldsymbol{v}(\boldsymbol{x}) \cdot 
    \boldsymbol{n}(\boldsymbol{x}) = v_{n}
    (\boldsymbol{x}) \quad \mbox{on} \; 
    \Gamma^{v} \\
    \label{Eqn:MDarcy_Pressure_BC}
    &p(\boldsymbol{x}) = p_0(\boldsymbol{x}) 
    \quad \mbox{on} \; \Gamma^{p} 
  \end{align}
\end{subequations}
where $\alpha(p)$ is the drag coefficient (which 
has dimensions of $\left[\mathrm{M} \mathrm{L}^{-3} 
\mathrm{T}^{-1}\right]$), $p_0(\boldsymbol{x})$ is 
the prescribed pressure, $v_{n}(\boldsymbol{x})$ 
is the prescribed normal component of the velocity, 
$\rho(\boldsymbol{x})$ denotes the density of the 
fluid, $\boldsymbol{b}(\boldsymbol{x})$ is the specific 
body force, and $\boldsymbol{n}(\boldsymbol{x})$ is the 
unit outward normal vector to $\partial \Omega$. The 
drag function is the ratio between viscosity of the 
fluid and the permeability. Herein we consider the 
following two forms for drag function
\begin{subequations}
  \label{Eqn:StabMDarcy_specific_alpha}
  \begin{align}
    \alpha(p) = \alpha_{0} (1 + \beta p) \\
    \alpha(p) = \alpha_{0} \exp[\beta p]
  \end{align}
\end{subequations}
The first equation can be considered as a two-term 
Taylor's series approximation of the second equation 
(which is based on the Barus' formula). It is, in 
general, not possible to obtain analytical solutions 
for the system of equations \eqref{Eqn:MDarcy_modified_Darcy}, 
especially for complex geometries. Hence, one may have to 
resort to numerical solutions. The main aim of this paper 
is to present a stabilized mixed formulation to solve the 
boundary value problem given by equations 
\eqref{Eqn:MDarcy_Equilibrium}--\eqref{Eqn:MDarcy_Pressure_BC}. 

Developing numerical formulations for the above equation falls in the realm of mixed methods \cite{Brezzi_Fortin}. 
As mentioned earlier, it is generally agreed upon that care should be taken in developing numerical formulations 
to avoid numerical instabilities. For example, under the classical Galerkin formulation (which is sometimes referred 
to as the classical mixed formulation) equal-order interpolation for velocity and pressure is not stable. A mathematical 
theory that addresses the stability issues with mixed methods is the Ladyzhenskaya-Babu\v ska-Brezzi (LBB) inf-sup 
condition. Specialized elements and appropriate function spaces have been developed for mixed problems (like 
Darcy equation, Stokes equation, incompressible Navier-Stokes) to either satisfy or circumvent the LBB stability 
condition. 

In the next subsection we present a new stabilized mixed formulation that is inherently more 
stable than the classical mixed formulation and can accommodate a larger variety of combinations 
of interpolation functions for velocity and pressure including the equal-order interpolation. Although 
the proposed formulation can accommodate  many combinations of interpolations for the velocity 
and pressure, in Section \ref{Sec:S4_StabMDarcy_Numerical_Results} we illustrate the performance 
of the proposed formulation using equal-order interpolation as other combinations are not computationally 
attractive.

\subsection{A new stabilized mixed formulation for modified Darcy equation}
\label{Subsec:S3_StabMDarcy_Stabilized_Formulations}
Let $\boldsymbol{w}(\boldsymbol{x})$ and $q(\boldsymbol{x})$ be test 
functions corresponding to the velocity and pressure, respectively. 
Let us define the following function spaces, which will be used in 
the remainder of this paper:
\begin{subequations}
  \begin{align}
    &\mathcal{V} := \left\{\boldsymbol{v}(\boldsymbol{x}) \in 
      \left(L_{2}(\Omega)\right)^{nd} \; \big| \; \mathrm{div}
      [\boldsymbol{v}] \in L_{2}(\Omega), \boldsymbol{v}
      (\boldsymbol{x}) \cdot \boldsymbol{n}(\boldsymbol{x}) 
      = v_n(\boldsymbol{x}) \; \mathrm{on} \; \Gamma^{v} \right\} \\
    &\mathcal{W} := \left\{\boldsymbol{w}(\boldsymbol{x}) \in 
      \left(L_{2}(\Omega)\right)^{nd} \; \big| \; \mathrm{div}
      [\boldsymbol{w}] \in L_{2}(\Omega), \boldsymbol{w}
      (\boldsymbol{x}) \cdot \boldsymbol{n}(\boldsymbol{x}) 
      = 0 \; \mathrm{on} \; \Gamma^{v} \right\} \\
    \label{Eqn:StabMDarcy_function_space_Q}
    &\mathcal{Q} \equiv H^{1}(\Omega)  
  \end{align}
\end{subequations}
In using equation \eqref{Eqn:StabMDarcy_function_space_Q} 
one assumes that $\Gamma^{p} \neq \emptyset$. If $\Gamma^{p} 
= \emptyset$ then, for well-posedness, $\mathcal{Q}$ should 
be taken as follows:
\begin{align*}
  \mathcal{Q} := \left\{p(\boldsymbol{x}) \in  H^{1}(\Omega) 
    \; \big| \; \int_{\Omega} p(\boldsymbol{x}) \; \mathrm{d} 
    \Omega = 0 \right\} 
\end{align*}
In addition, for the case of $\Gamma^{p} = \emptyset$, we 
need to satisfy the following compatibility condition on 
the prescribed data: 
\begin{align}
\int_{\partial \Omega} v_n(\boldsymbol{x}) \; \mathrm{d} \Gamma = 0 
\end{align}
For convenience, let us denote the standard $L_{2}$ inner product 
defined over spatial domain $K$ as $(\cdot;\cdot)_K$. That is,
\begin{align}
  (a;b)_{K} \equiv \int_{K} a \cdot b \; \mathrm{d}K \quad \forall a, b 
\end{align}
For simplicity, the subscript $K$ will be dropped if $K = \Omega$. 

In References \cite{Masud_Hughes_CMAME_2002_v191_p4341,
Nakshatrala_Turner_Hjelmstad_Masud_CMAME_2006_v195_p4036} 
a stabilized formulation based on the variational multiscale 
(VMS) formalism has been proposed and analyzed for the case 
of Darcy equation (which assumes constant drag coefficient). 
In this paper we extend the variational multiscale formulation 
to the modified Darcy equation (which gives rise to nonlinear 
equations). We shall linearize the governing equations 
\eqref{Eqn:MDarcy_Equilibrium}--\eqref{Eqn:MDarcy_Pressure_BC} 
using a fixed-point procedure as follows: 
\begin{subequations}
  \begin{align}
    &\alpha(p^{(i-1)}) \boldsymbol{v}^{(i)} + \mathrm{grad}[p^{(i)}] 
    = \rho(\boldsymbol{x}) 
    \boldsymbol{b} (\boldsymbol{x}) \quad \mathrm{in} \; \Omega \\
    &\mathrm{div}[\boldsymbol{v}^{(i)}] = 0 \quad \mathrm{in} \; \Omega \\
    &\boldsymbol{v}^{(i)}(\boldsymbol{x}) \cdot \boldsymbol{n}(\boldsymbol{x}) = 
    v_n(\boldsymbol{x}) \quad \mathrm{in} \; \Gamma^{v} \\
    &p^{(i)}(\boldsymbol{x}) = p_{0}(\boldsymbol{x}) \quad \; \mathrm{on} \; \Gamma^{p}
  \end{align}
\end{subequations}  
Based on the derivation presented in Reference 
\cite{Nakshatrala_Turner_Hjelmstad_Masud_CMAME_2006_v195_p4036}, 
a stabilized mixed formulation for the above equations based on 
the VMS formalism can be written as follows: Find $\boldsymbol{v}^{(i)}
(\boldsymbol{x}) \in \mathcal{V}$ and $p^{(i)}(\boldsymbol{x}) 
\in \mathcal{Q}$ such that we have
\begin{align}
  \label{Eqn:StabMDarcy_VMS1_formulation_step}
  \mathcal{G}_{\mathrm{stab}}(\boldsymbol{w},q;\boldsymbol{v}^{(i)},p^{(i)};p^{(i-1)}) 
  = \mathcal{L}_{\mathrm{stab}}(\boldsymbol{w},q;p^{(i-1)}) \quad \forall \boldsymbol{w}(\boldsymbol{x}) 
  \in \mathcal{W}, q(\boldsymbol{x}) \in \mathcal{Q}
\end{align}
where the functionals $\mathcal{G}_{\mathrm{stab}}$ and 
$\mathcal{L}_{\mathrm{stab}}$ are, respectively, defined 
as follows:
\begin{subequations}
  \begin{align}
    \label{Eqn:StabMDarcy_Gstab}
    \mathcal{G}_{\mathrm{stab}}(\boldsymbol{w},q;\boldsymbol{v},p;\tilde{p}) 
    &:= (\boldsymbol{w}; \alpha(\tilde{p}) \boldsymbol{v}) - (\mathrm{div}
    [\boldsymbol{w}];p) - (q;\mathrm{div}[\boldsymbol{v}]) \notag \\ 
    &- \frac{1}{2} \left(\alpha(\tilde{p}) \boldsymbol{w} + \mathrm{grad}[q]; 
      \alpha^{-1}(\tilde{p}) \left(\alpha(\tilde{p}) \boldsymbol{v} + 
        \mathrm{grad}[p] \right) \right) \\
    \label{Eqn:StabMDarcy_Lstab}
    \mathcal{L}_{\mathrm{stab}}(\boldsymbol{w},q;\tilde{p}) &:= (\boldsymbol{w};
    \rho \boldsymbol{b}) - (\boldsymbol{w} \cdot \boldsymbol{n}; p_0)_{\Gamma^{p}} 
    - \frac{1}{2} \left(\alpha(\tilde{p}) \boldsymbol{w} + \mathrm{grad}[q]; 
      \alpha^{-1}(\tilde{p}) \rho \boldsymbol{b} \right)
  \end{align}
\end{subequations}
The terms in the second line of equation \eqref{Eqn:StabMDarcy_Gstab} 
and the last term in \eqref{Eqn:StabMDarcy_Lstab} are referred to 
as stabilization terms. The factor $1/2$ (in front of these terms) 
is the stabilization parameter. The last term in equation 
\eqref{Eqn:StabMDarcy_Lstab} is the stabilization term due 
to the body force. The stability of the above formulation can 
be inferred using the mathematical proof outlined in Reference 
\cite{Masud_Hughes_CMAME_2002_v191_p4341}.

The proposed numerical algorithm for solving the original system 
of equations \eqref{Eqn:MDarcy_modified_Darcy} is given in Algorithm 
\ref{Algo:StabMDarcy_VMS1_formulation}, and will be based on the 
above stabilized formulation. Since the formulation is based on 
the fixed-point linearization, the rate of convergence of the 
algorithm (if it converges) will be linear with respect to 
iteration number. Moreover, when $\alpha$ is a constant (that 
is, for classical Darcy equation), the proposed formulation will 
converge in one iteration. In a finite element implementation, 
the norm in the stopping criterion (see line 10 in Algorithm 
\ref{Algo:StabMDarcy_VMS1_formulation}) can be taken as the 
Euclidean norm of the vector containing nodal pressures. 

\begin{algorithm}
  \caption{\textsf{A new stabilized formulation for modified Darcy equation}}
  \label{Algo:StabMDarcy_VMS1_formulation}
  \begin{algorithmic}[1]
    \STATE Input: $\epsilon_{\mathrm{TOL}}$, MAXITERS
    \STATE Output: $\boldsymbol{v}_h(\boldsymbol{x})$ and $p_h(\boldsymbol{x})$
    \STATE Guess $p_h^{(0)}(\boldsymbol{x}) \in \mathcal{Q}_h$
    \STATE Initialize: $i \leftarrow 1$
    \WHILE {true}
    \IF{$i > \mathrm{MAXITERS}$} 
    \PRINT Iterative scheme did not converge. RETURN 
    \ENDIF
    \STATE Find $\boldsymbol{v}_h^{(i)}(\boldsymbol{x}) \in \mathcal{V}_{h}$ 
    and $p_h^{(i)}(\boldsymbol{x}) \in \mathcal{Q}_h$ such that 
    \begin{align*}
      \mathcal{G}_{\mathrm{stab}}(\boldsymbol{w}_h,q_h; \boldsymbol{v}_h^{(i)},p_h^{(i)};
      p_h^{(i-1)}) = \mathcal{L}_{\mathrm{stab}}(\boldsymbol{w}_h,q_h;p_h^{(i-1)}) \quad 
      \forall \boldsymbol{w}_h(\boldsymbol{x}) \in \mathcal{W}_h, \; q_h(\boldsymbol{x}) 
      \in \mathcal{Q}_h
    \end{align*}
    \COMMENT{Check for convergence based on the given tolerance}
    \IF{$\|p_h^{(i)} - p_h^{(i-1)} \| < \epsilon_{\mathrm{TOL}}$}
    \PRINT Iterative scheme converged. 
    \STATE $\boldsymbol{v}_h(\boldsymbol{x}) \leftarrow 
    \boldsymbol{v}_h^{(i)}(\boldsymbol{x})$ and $p_h(\boldsymbol{x}) 
    \leftarrow p_h^{(i)}(\boldsymbol{x})$ 
    \STATE RETURN
    \ENDIF
    \STATE $i \leftarrow i + 1$
    \ENDWHILE
  \end{algorithmic}
\end{algorithm}

%
\section{REPRESENTATIVE NUMERICAL RESULTS}
\label{Sec:S4_StabMDarcy_Numerical_Results}
In this section, we present a wide variety of representative 
numerical results to illustrate the performance of the proposed 
stabilized formulation. The first few test problems are used to 
illustrate the accuracy of the proposed formulation on various 
computational grids as for these test problems analytical solutions 
can be easily obtained. One test problem is used to illustrate 
the robustness of the formulation in the case of heterogeneous 
medium. To the end of the section, the performance of the proposed 
formulation in a parallel setting is illustrated using a large-scale 
simulation related to carbon-dioxide sequestration.
\emph{It should be noted that in all our numerical simulations 
we have employed low-order Lagrange finite element approximations 
(i.e., $p = 1$), and equal-order interpolation for the velocity 
and pressure.}
 
We shall first non-dimensionalize the governing equations by 
employing the same non-dimensionalization procedure as given 
in Reference \cite{Nakshatrala_Rajagopal_IJNMF_arxiv_2010}. 
We shall present the details for completeness. We define 
the following non-dimensional quantities (which have a 
superposed bar):
\begin{align}
  \bar{\boldsymbol{x}} = \frac{\boldsymbol{x}}{L}, \; 
  \bar{\boldsymbol{v}} = \frac{\boldsymbol{v}}{V}, \; 
  \bar{p} = \frac{p}{P}, \; 
  \bar{\alpha} = \frac{\alpha}{\alpha_{\mathrm{ref}}}, \; 
  \bar{\alpha}_0 = \frac{\alpha_0}{\alpha_{\mathrm{ref}}},\; 
  \bar{\rho} = \frac{\rho}{\rho_{\mathrm{ref}}}, \; 
  \bar{\boldsymbol{b}} = \frac{\boldsymbol{b}}{B}, \; 
  \bar{\beta} = \beta P, \; 
  \bar{v}_n = \frac{v_n}{V}, \; 
  \bar{p}_0 = \frac{p_0}{P}
\end{align}
where $L$, $V$, $P$, $\alpha_{\mathrm{ref}}$, $\rho_{\mathrm{ref}}$ 
and $B$ respectively denote reference length, velocity, pressure, 
drag coefficient, density and specific body force. The gradient 
and divergence operators with respect to $\bar{\boldsymbol{x}}$ 
are denoted as ``$\overline{\mathrm{grad}}$" and 
``$\overline{\mathrm{div}}$", respectively. The scaled domain 
$\Omega_{\mathrm{scaled}}$ is defined as follows: a point 
in space with position vector $\bar{\boldsymbol{x}} \in 
\Omega_{\mathrm{scaled}}$ corresponds to the same point with 
position vector given by $\boldsymbol{x} = \bar{\boldsymbol{x}} 
L \in \Omega$ . Similarly, one can define the scaled boundaries: 
$\partial \Omega_{\mathrm{scaled}}$, $\Gamma^{\mathrm{D}}_{\mathrm{scaled}}$, 
and $\Gamma^{\mathrm{N}}_{\mathrm{scaled}}$. The above non-dimensionalization 
gives rise to two dimensionless parameters
\begin{align}
  \mathcal{A} := \frac{\alpha_{\mathrm{ref}} V L}{P} \quad \mathrm{and} \quad \mathcal{C}:= \frac{\rho_{\mathrm{ref}} L B}{P}
\end{align}
A corresponding non-dimensionalized form of the drag functions given in equation 
\eqref{Eqn:StabMDarcy_specific_alpha} can be written as
\begin{subequations} 
  \begin{align}
    \label{Eqn:StabMDarcy_non_dimensional_linear}
    \bar{\alpha}(\bar{p}) &= \bar{\alpha}_0 
    \left(1 + \bar{\beta} \bar{p}\right) \\
    \label{Eqn:StabMDarcy_non_dimensional_Barus}
    \bar{\alpha}(\bar{p}) &= \bar{\alpha}_0 
    \exp[\bar{\beta} \bar{p}] 
  \end{align}
\end{subequations}
We shall write a non-dimensional form of the modified Darcy equation as 
\begin{subequations}
  \begin{align}
    &\mathcal{A} \; \bar{\alpha}(\bar{p}) \bar{\boldsymbol{v}} + \overline{\mathrm{grad}} [\bar{p}] = 
    \mathcal{C} \; \bar{\rho} \; \bar{\boldsymbol{b}}(\bar{\boldsymbol{x}}) \quad \; \mathrm{in} \; 
    {\Omega}_{\mathrm{scaled}} \\
    &\overline{\mathrm{div}}[\bar{\boldsymbol{v}}] = 0 \quad \; \mathrm{in} \; {\Omega}_{\mathrm{scaled}} \\
    &\bar{\boldsymbol{v}}(\bar{\boldsymbol{x}}) \cdot \bar{\boldsymbol{n}}(\bar{\boldsymbol{x}}) = 
    \bar{v}_n(\bar{\boldsymbol{x}}) \quad \; \mathrm{on} \; {\Gamma}^{\mathrm{N}}_{\mathrm{scaled}} \\
    &\bar{p}(\bar{\boldsymbol{x}}) = \bar{p}_0 (\bar{\boldsymbol{x}}) \quad \; \mathrm{on} \; 
    {\Gamma}^{\mathrm{D}}_{\mathrm{scaled}}
  \end{align}
\end{subequations}
where $\bar{\boldsymbol{n}}(\bar{\boldsymbol{x}})$ is the unit 
outward normal to the boundary $\partial {\Omega}_{\mathrm{scaled}}$. 

\subsection{One-dimensional problem}
\begin{figure}
  \psfrag{p1}{$p(x = 0) = p_1$}
  \psfrag{p2}{$p(x = 1) = p_2$}
  \psfrag{1}{$1.0$}
  \psfrag{x}{$x$}
  \includegraphics[scale=0.8]{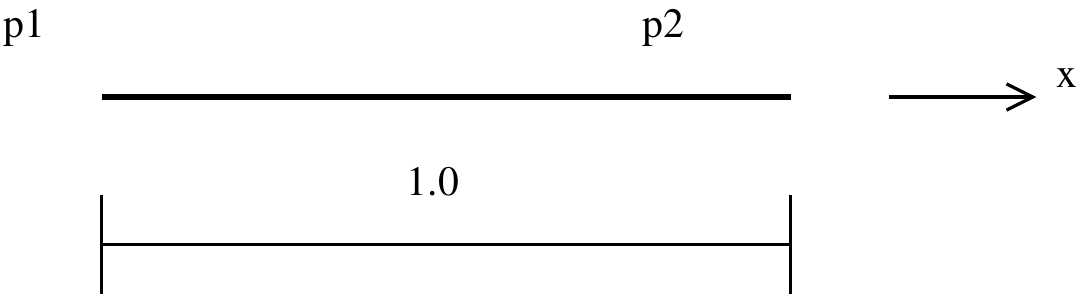}
  \caption{A pictorial description of the one-dimensional problem.\label{Fig:StabMDarcy_OneD_problem_description}}
\end{figure}

\begin{figure}
\centering
\subfigure{
\includegraphics[scale=0.35]{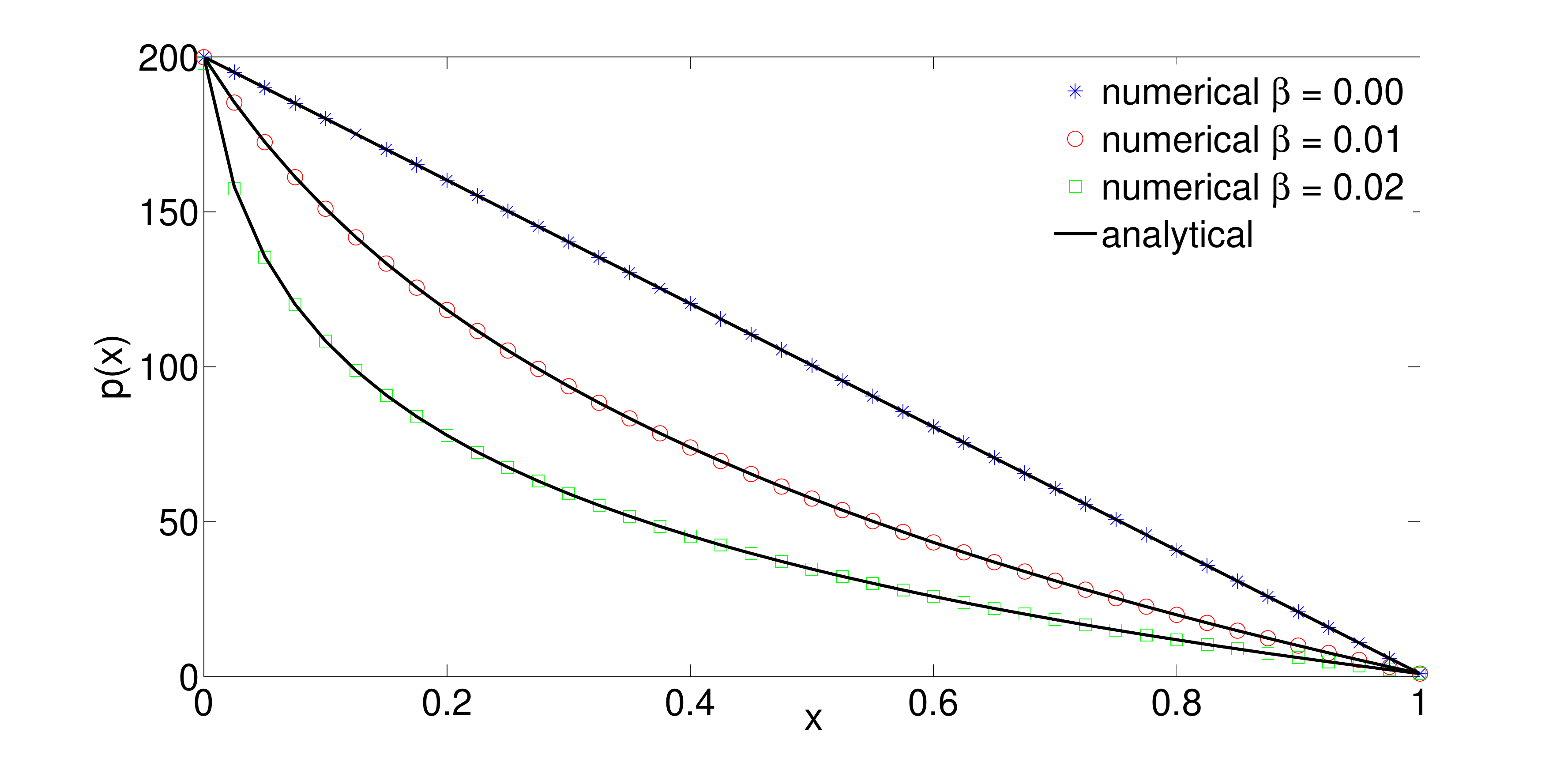}}
\subfigure{
\includegraphics[scale=0.35]{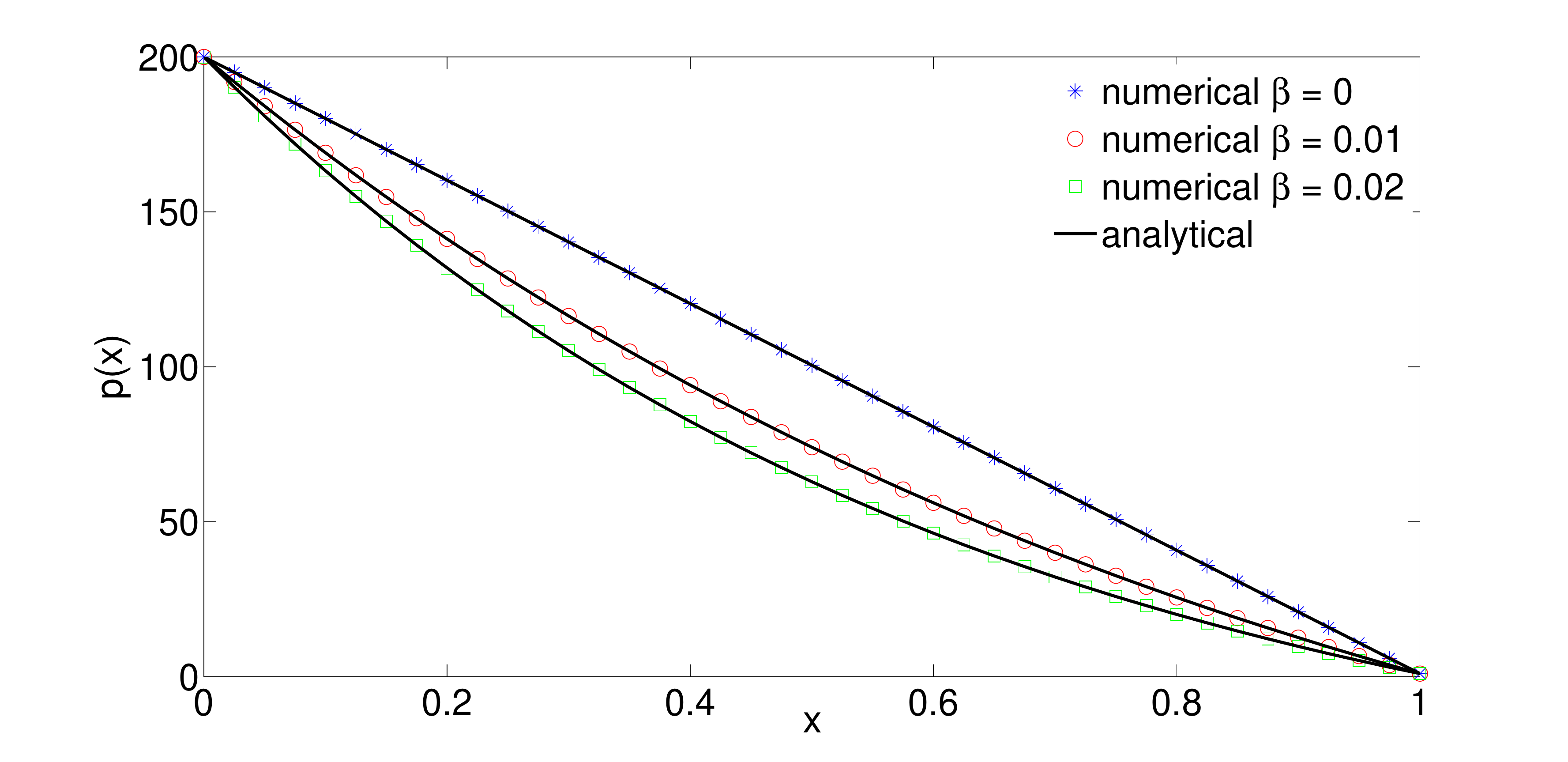}}
\caption{One-dimensional problem: In the top figure Barus formula is employed, and 
in the bottom figure linear variation of viscosity with respect to pressure is employed. 
Pressure is plotted along $x$ for various values of $\beta$. In this numerical example 
we have taken $\bar{\alpha}_0 = 1$, $\bar{p}(\bar{x} = 0) = 200$, $\bar{p}(\bar{x} = 1) = 1$, 
and $\epsilon_{\mathrm{TOL} } = 10^{-10}$. 
\label{Fig:StabMDarcy_OneD_Pressure}}
\end{figure}

\begin{figure}
  \centering
  \includegraphics[scale=0.35]{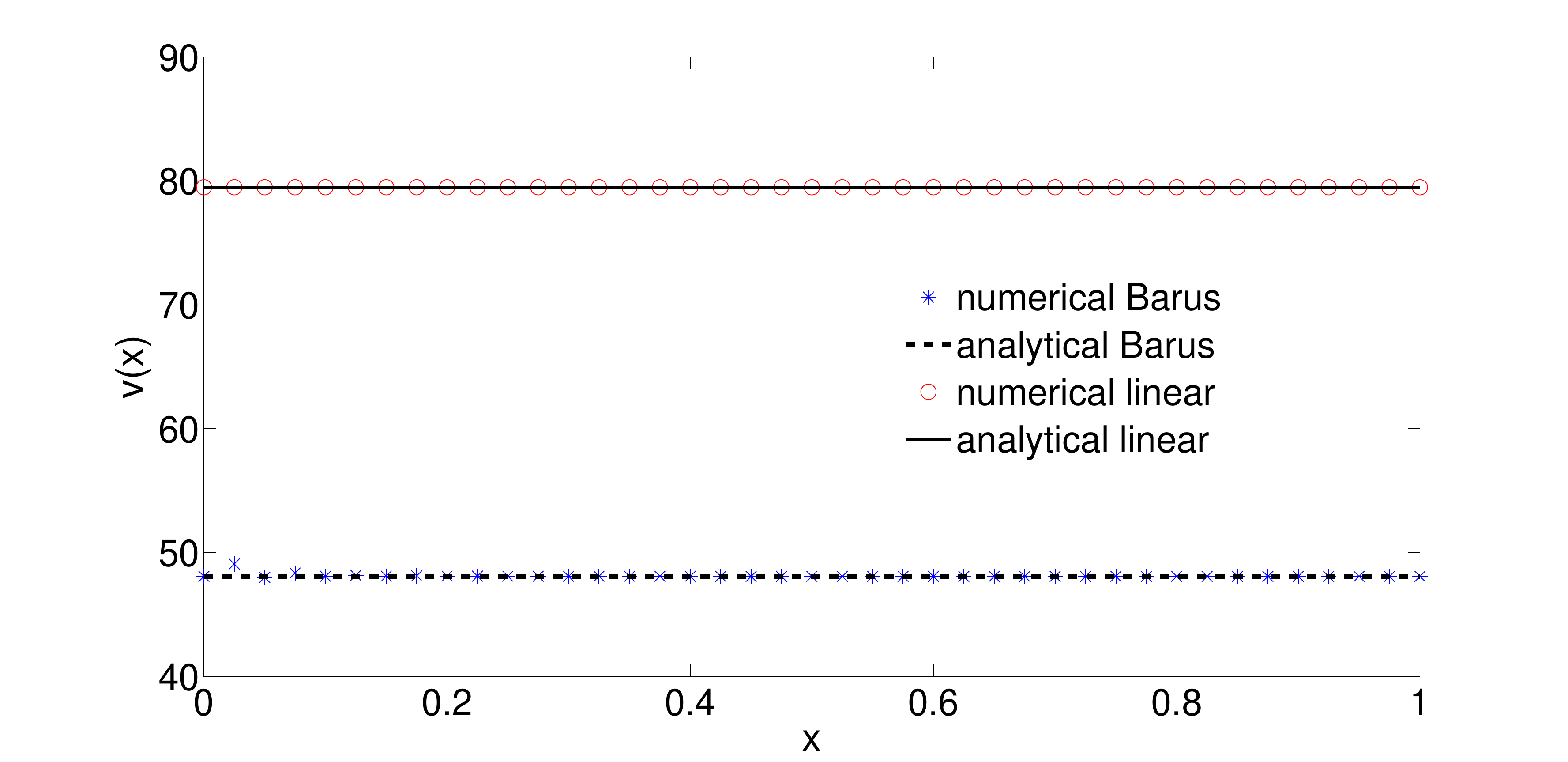}
\caption{One-dimensional problem: Velocity is plotted along $x$ 
  for various values of $\beta$. In this numerical example 
  we have taken $\bar{\alpha}_0 = 1$, $\bar{\beta} = 0.02$, 
  $\bar{p}(\bar{x} = 0) = 200$, $\bar{p}(\bar{x} = 1) = 1$, 
  and $\epsilon_{\mathrm{TOL} } = 10^{-10}$. 
  \label{Fig:StabMDarcy_OneD_Velocity}}
\end{figure}

\begin{figure}
  \label{Fig:StabMDarcy_OneD_iters_vs_beta}
  \psfrag{beta}{$\beta$}
  \includegraphics[scale=0.37,clip]{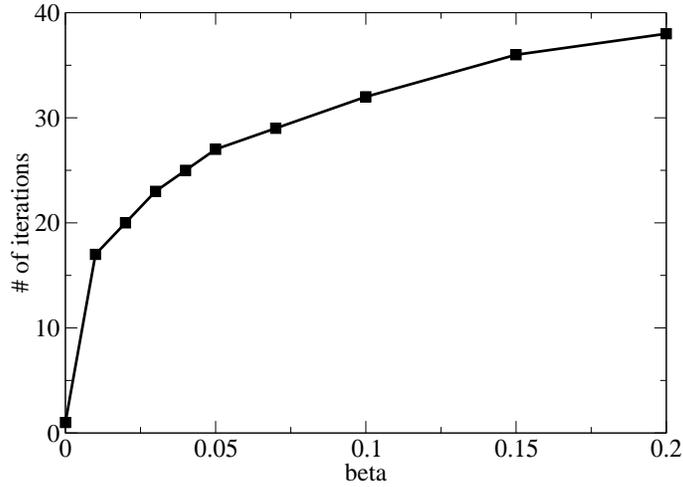}
  \caption{One-dimensional problem: This figure shows the number 
    of iterations taken by the numerical algorithm for various values 
    of $\beta$, which is an indicator of the strength of the nonlinearity.}
\end{figure}

\begin{figure}
  \label{Fig:StabMDarcy_OneD_res_vs_iters}
  \includegraphics[scale=0.5]{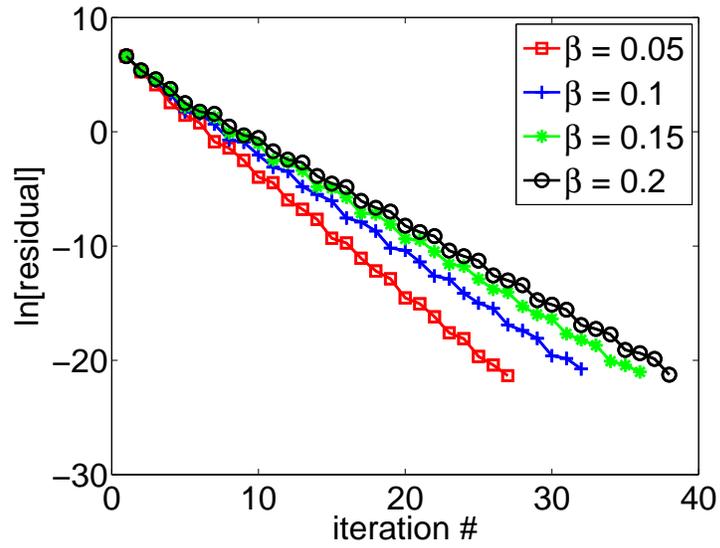}
  \caption{One-dimensional problem: This figure shows the 
  variation of the residual with respect to iteration number 
  for various values of $\beta$. Note that the $y$-axis is 
  natural logarithm of the residual.}
\end{figure}

We shall test the proposed stabilized mixed formulation using 
a simple one-dimensional problem. Consider the computational 
domain to be of unit length, and pressures of $\bar{p}_1$ and 
$\bar{p}_2$ are respectively prescribed at the left and right 
ends of the unit domain (see Figure 
\ref{Fig:StabMDarcy_OneD_problem_description}). We neglect the 
body force. The governing equations for this test problem can 
be written as 
\begin{subequations}
  \begin{align}
    &\mathcal{A} \bar{\alpha}(\bar{p}) \bar{v}(\bar{x}) + \frac{\mathrm{d} 
      \bar{p}}{\mathrm{d} \bar{x}} = 0 \quad \mathrm{in} \; (0,1) \\
    &\frac{\mathrm{d} \bar{v}}{d \bar{x}} = 0 \quad \mathrm{in} \; (0,1) \\
    &\bar{p}(\bar{x} = 0) = \bar{p}_1, \quad \bar{p}(\bar{x} = 1) = \bar{p}_2
  \end{align}
\end{subequations}
For this test problem, the analytical solutions with the drag 
functions defined in equation \eqref{Eqn:StabMDarcy_specific_alpha} 
can be written as 
\begin{subequations}
  \label{Eqn:StabMDarcy_OneD_solutions}
  \begin{align}
    \mbox{for the case} \; \bar{\alpha}(\bar{p}) = \bar{\alpha}_0 (1 + \bar{\beta} \bar{p}) \; \left\{
      \begin{array}{l}
        \bar{p}(\bar{x}) = \frac{1}{\bar{\beta}} \left[\left(1 + \bar{\beta}\bar{p}_1\right)^{1 - \bar{x}} 
          \left(1 + \bar{\beta} \bar{p}_2\right)^{\bar{x}} - 1\right] \\
        \bar{v}(\bar{x}) = \frac{-1}{\mathcal{A} \bar{\alpha}_0 \bar{\beta}} \ln 
        \left[\frac{1 + \bar{\beta}\bar{p}_2}{1 + \bar{\beta}\bar{p}_1}\right]
      \end{array} \right.
  \end{align}
  \begin{align}
    \mbox{for the case} \; \bar{\alpha}(\bar{p}) = \bar{\alpha}_0 \exp[\bar{\beta} \bar{p}] \; \left\{
      \begin{array}{l}
        \bar{p}(\bar{x}) = \frac{-1}{\bar{\beta}} \ln \left\{ (1 - \bar{x}) \exp\left[-\bar{\beta}\bar{p}_1
          \right] + \bar{x} \exp\left[-\bar{\beta} \bar{p}_2\right] \right\} \\
        \bar{v}(\bar{x}) = \frac{1}{\mathcal{A} \bar{\alpha}_0 \bar{\beta}} \left\{
          \exp\left[-\bar{\beta} \bar{p}_2\right] - \exp\left[-\bar{\beta} \bar{p}_1\right] \right\}
      \end{array} \right.
  \end{align}
\end{subequations}
In Figures \ref{Fig:StabMDarcy_OneD_Pressure} and 
\ref{Fig:StabMDarcy_OneD_Velocity} we compare the 
numerical solutions against the analytical solutions 
for various drag functions. In all the cases considered, 
the proposed numerical formulation performed well. It is 
important to note that steep gradients in the pressure 
occur near the boundary for high values of $\beta$ 
(which is an indicator of the strength of nonlinearity).
Figure \ref{Fig:StabMDarcy_OneD_iters_vs_beta} shows the 
number of iterations taken by the numerical algorithm with 
respect to $\beta$. As expected, the number of iterations 
increases with respect to $\beta$. 
Figure \ref{Fig:StabMDarcy_OneD_res_vs_iters} shows the 
variation of the residual with respect to the iteration 
number for various value of $\beta$. The residual decreases 
monotonically with the iteration number.

\subsection{Two-dimensional constant flow problem}
\begin{figure}
  \centering
  \subfigure{
    \includegraphics[scale=0.35]{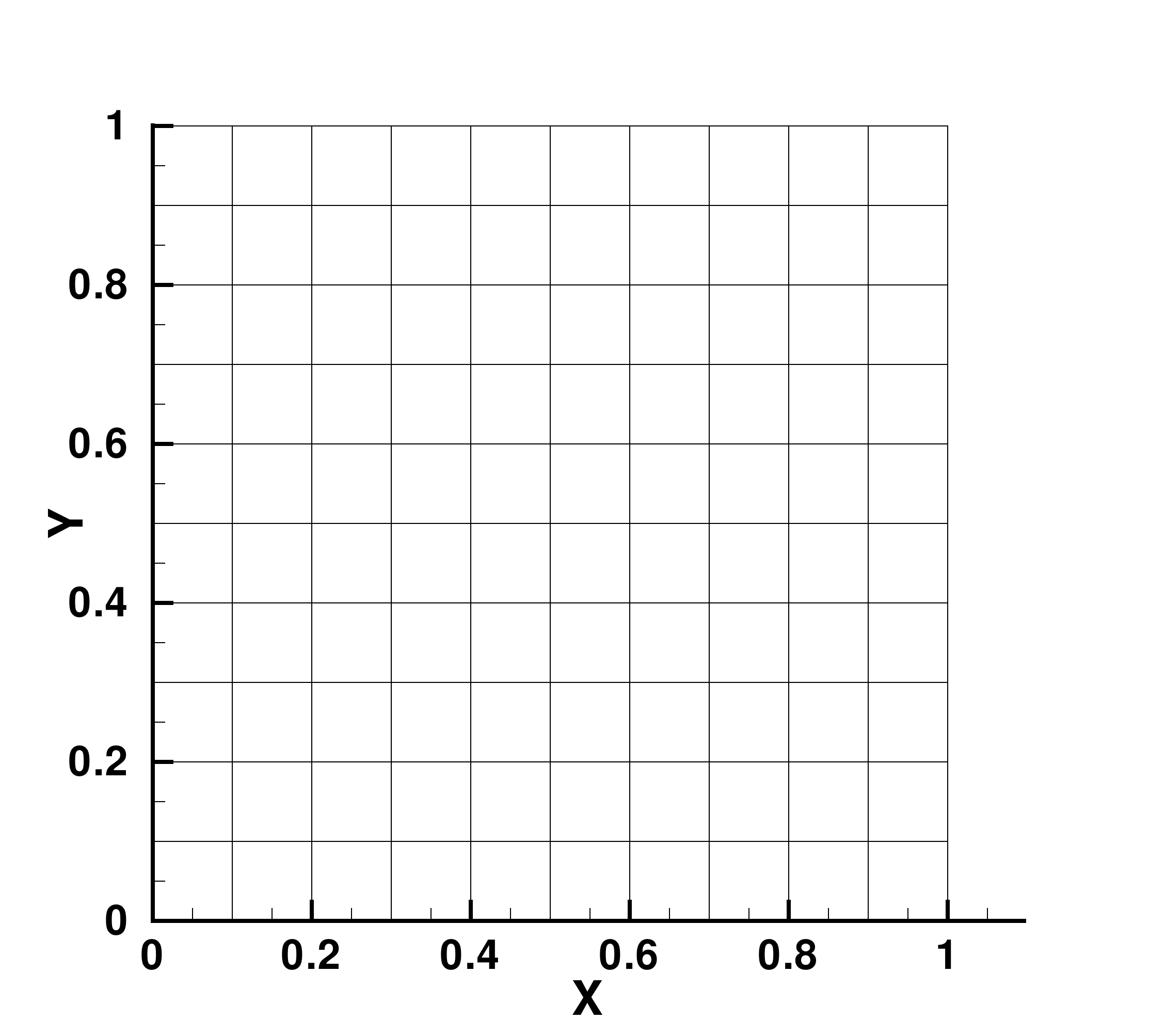}}
  \subfigure{
    \includegraphics[scale=0.35]{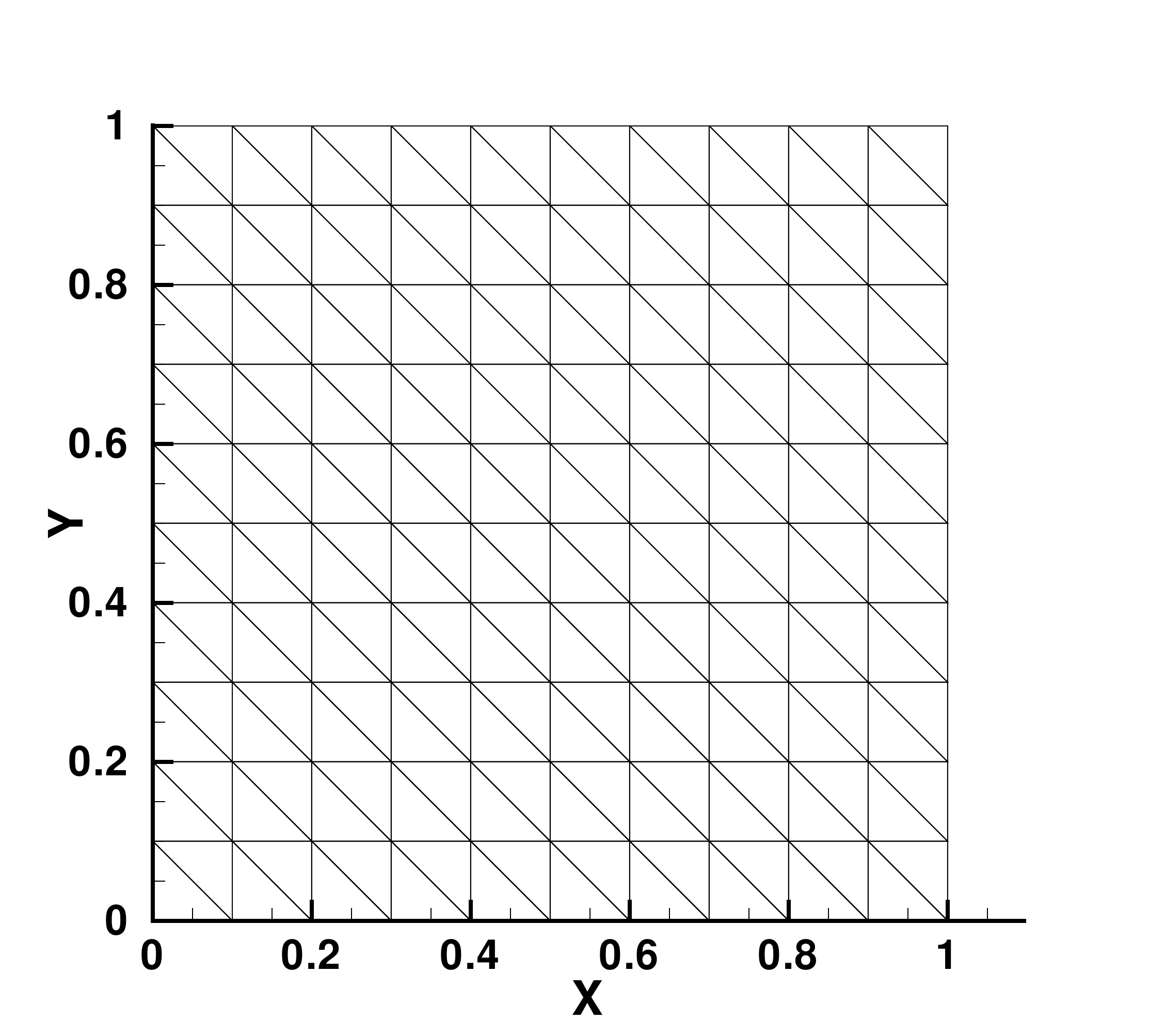}}
  \caption{Two-dimensional constant flow problem: Four-node quadrilateral 
    (left) and three-node triangular (right) meshes used in numerical 
    simulations. \label{Fig:StabMDarcy_TwoD_meshes}}
\end{figure}

\begin{figure}
\centering
\subfigure{
\includegraphics[scale=0.35]{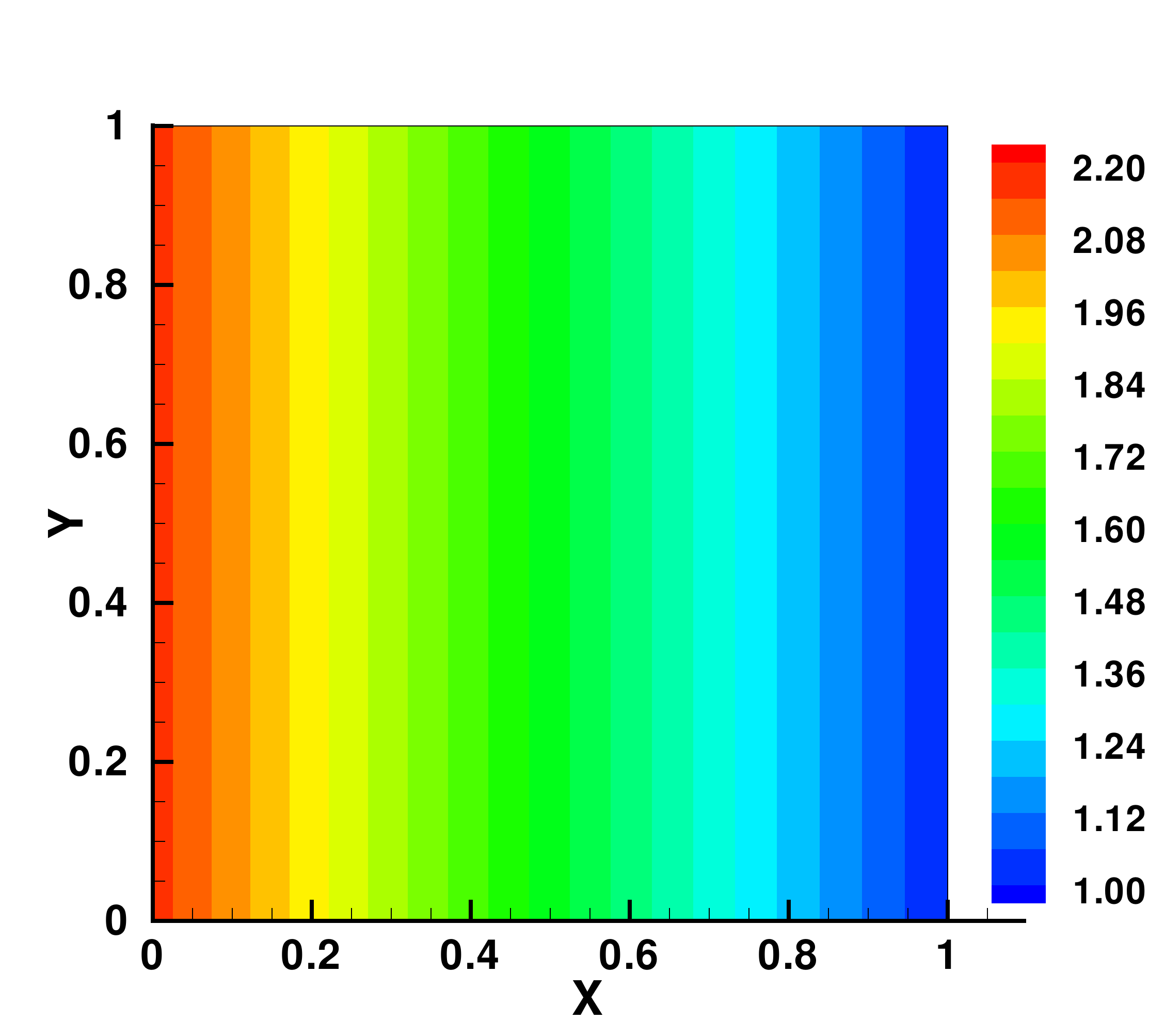}}
\subfigure{
\includegraphics[scale=0.35]{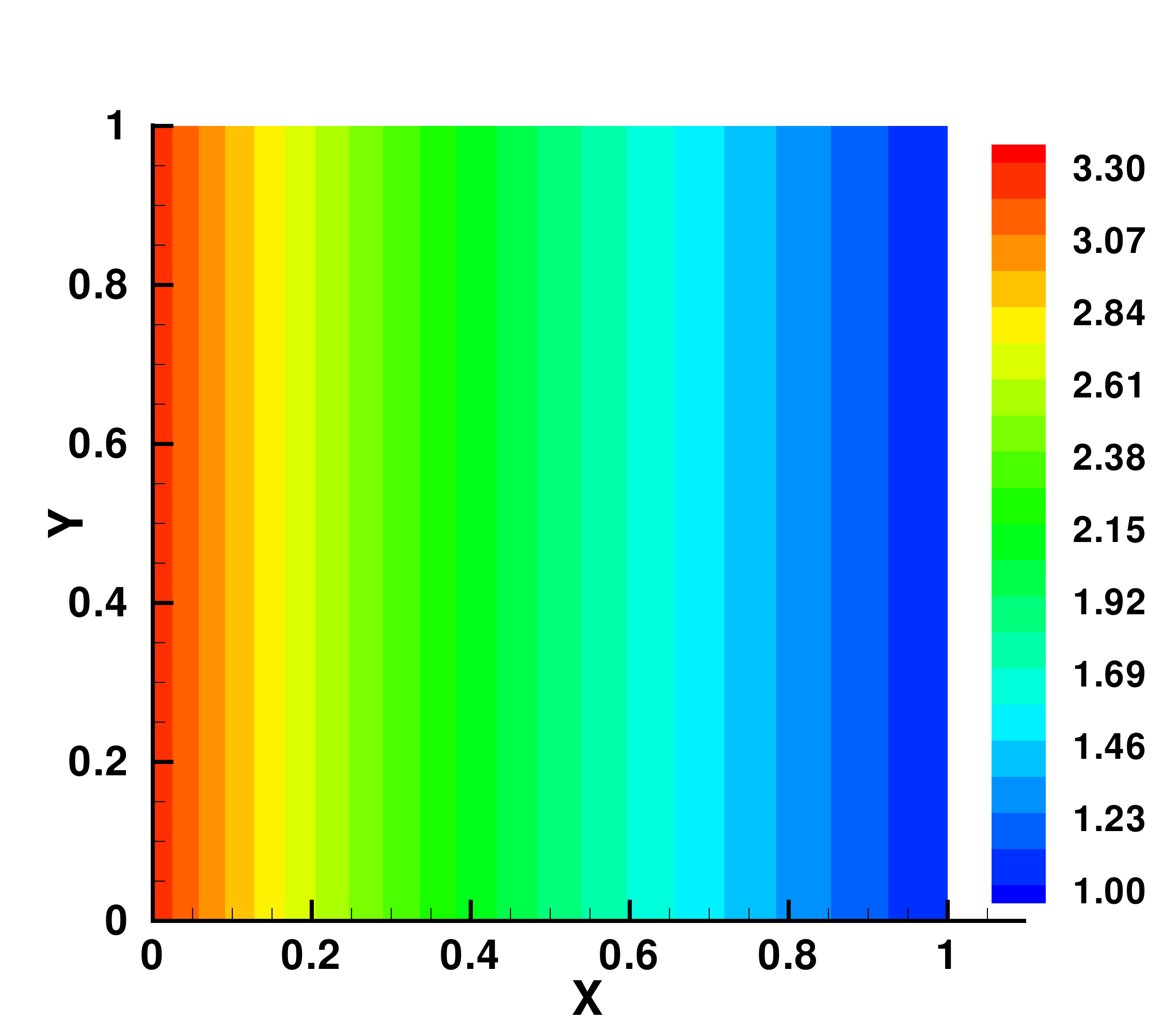}}
\subfigure{
\includegraphics[scale=0.35]{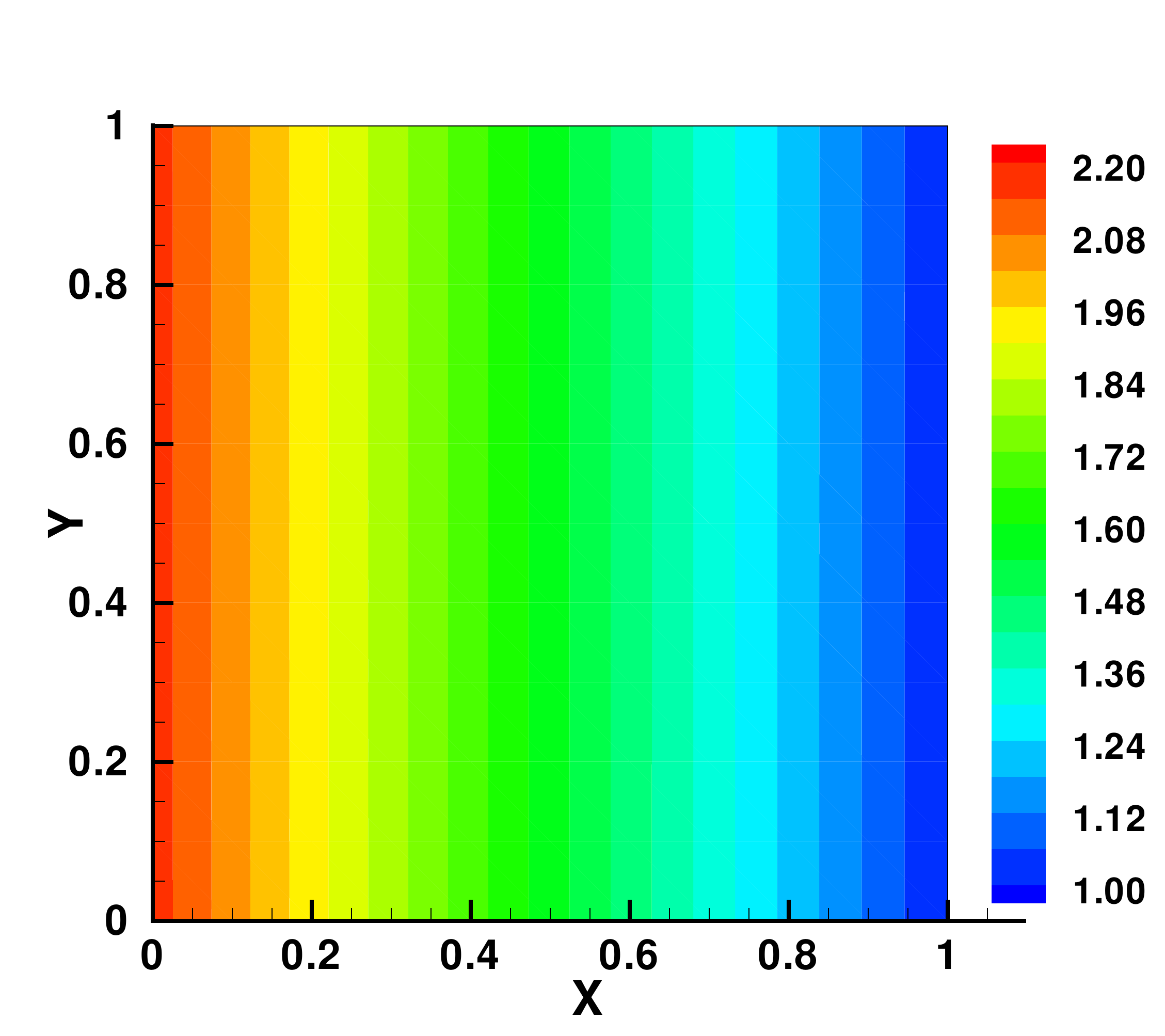}}
\subfigure{
\includegraphics[scale=0.35]{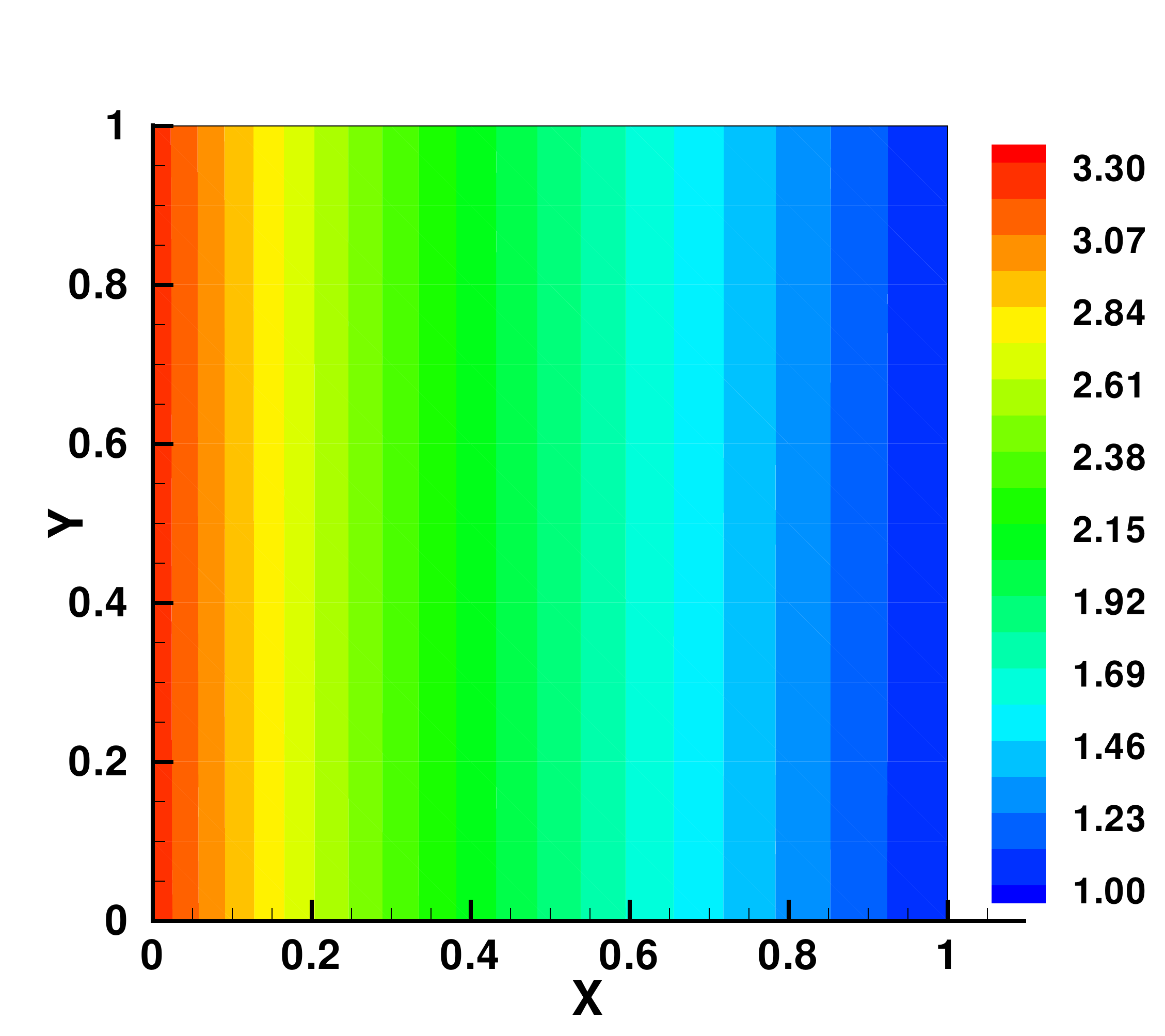}}
\caption{Two-dimensional constant flow problem: Barus formula is employed with $\alpha_0 = 1$, 
and the tolerance is taken to be $\epsilon_{\mathrm{TOL}} = 10^{-10}$. For the left figure we have 
used $\beta = 0.1$ and for the right figure we have used $\beta = 0.4$. Four-node quadrilateral mesh 
is used for top figures, and three-node triangular mesh is used for bottom figures. The numerical results 
matched well with the analytical solution. \label{Fig:StabMDarcy_TwoD_Pressure}}
\end{figure}

The computational domain is a bi-unit square domain, $\bar{\Omega} = [0, 1] \times [0, 1]$. Barus formula 
with $\bar{\alpha}_0 = 1$ is employed with two different values of $\bar{\beta} = 0.1$ and $\bar{\beta} = 0.4$. 
The body force is neglected. The left and right sides of the domain are prescribed with $\bar{v}_{\bar{x}} = 1$. 
The top and bottom sides of the domain are prescribed with $\bar{v}_{\bar{y}} = 0$. The analytical solution is 
given by 
\begin{subequations}
  \begin{align}
    &\bar{v}_{\bar{x}}(\bar{x},\bar{y}) = 1, \quad 
    \bar{v}_{\bar{y}}(\bar{x},\bar{y}) = 0 \\
    &\bar{p}(\bar{x},\bar{y}) = p_0 - \frac{1}{\bar{\beta}} 
    \ln\left[1 - \bar{\alpha}_0 \bar{\beta} (1 - \bar{x}) 
      \exp[\bar{\beta} p_0]\right]
  \end{align}
\end{subequations}
We prescribed $p_0 = 1$ at the top-right corner of the computational 
domain. The finite element meshes used in the numerical simulations 
are shown in Figure \ref{Fig:StabMDarcy_TwoD_meshes}. The pressure 
contours using the proposed stabilized formulation are shown in 
Figure \ref{Fig:StabMDarcy_TwoD_Pressure}.

\subsection{The quarter five-spot problem}
\begin{figure}
  \psfrag{vx}{$v_x = 0$}
  \psfrag{vy}{$v_y = 0$}
  \includegraphics[scale=0.75]{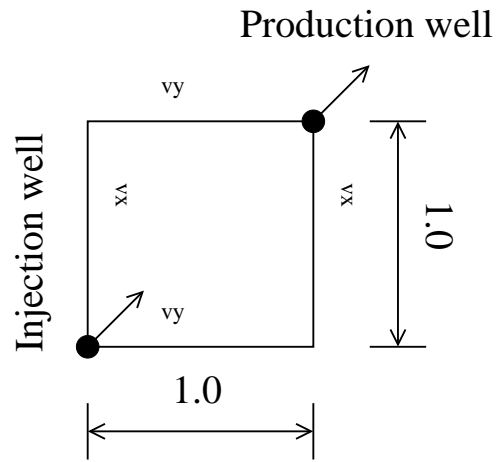}
  \caption{A pictorial description of the quarter five-spot problem.
    \label{Fig:StabMDarcy_Five_spot_description}}
\end{figure}

\begin{figure}
  \subfigure{
    \includegraphics[scale=0.5]{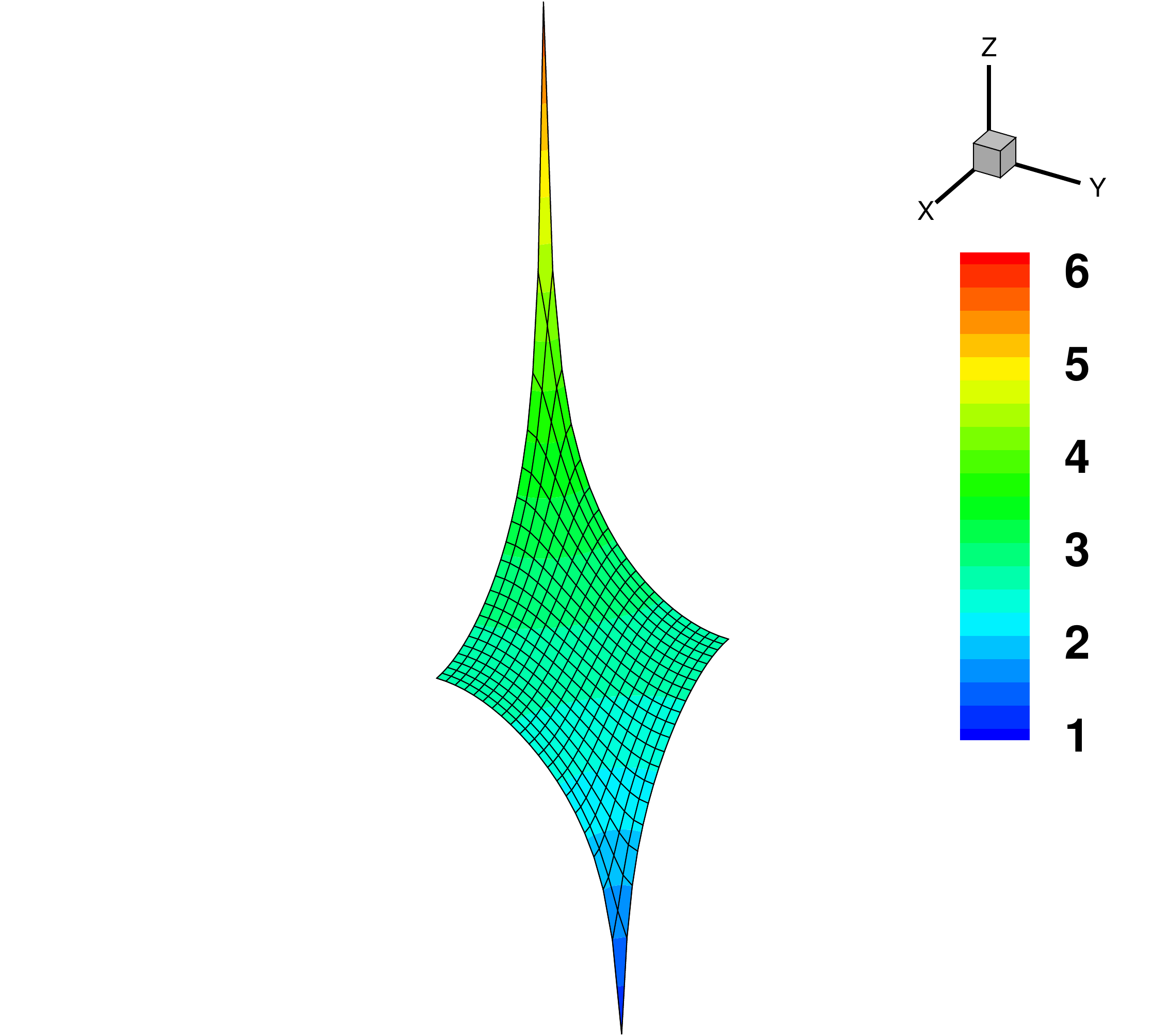}}
  \subfigure{
    \includegraphics[scale=0.5]{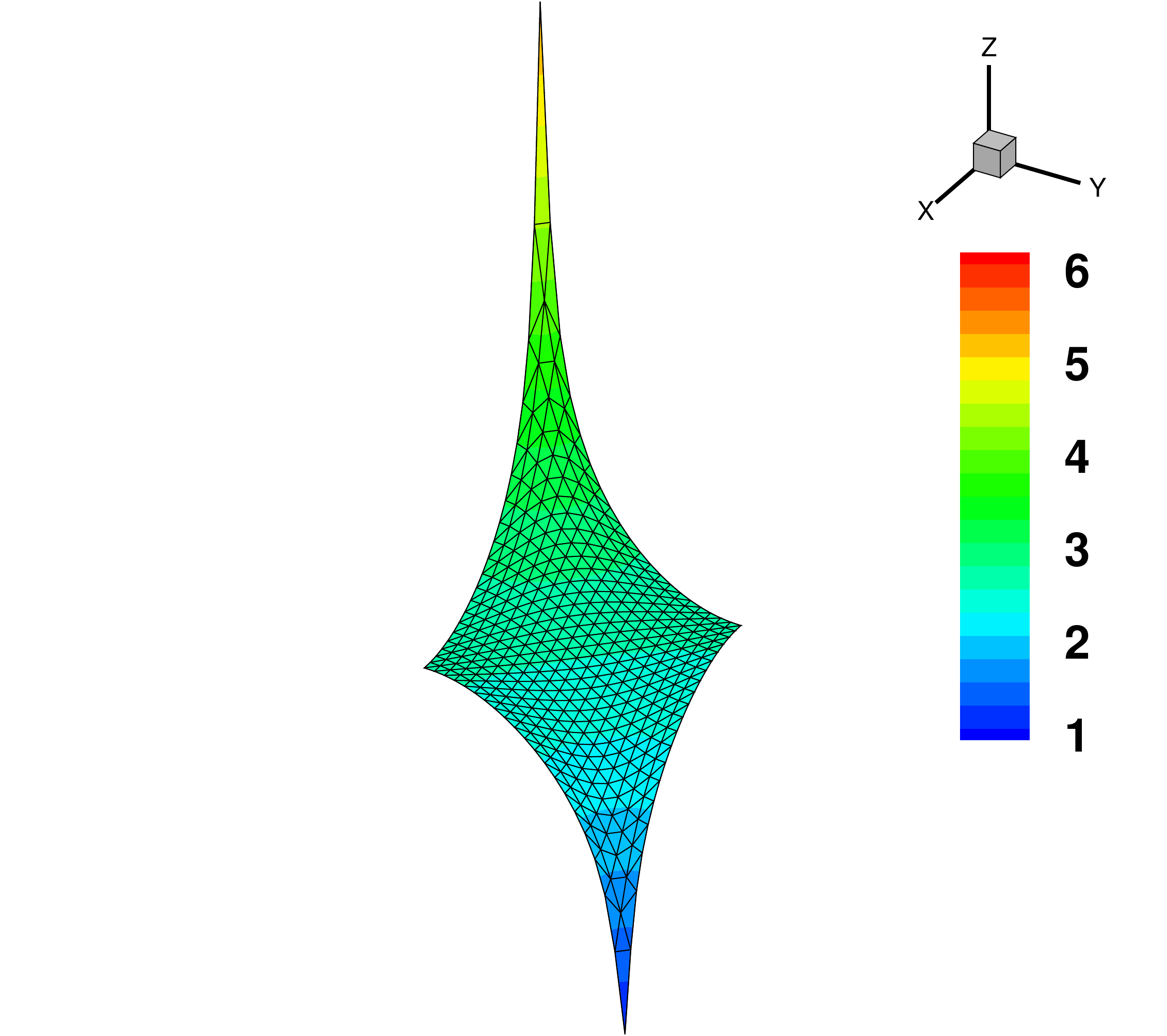}}
  \caption{Five-spot problem: Pressure contours using four-node quadrilateral (top) and three-node triangular 
  (bottom) elements. We have used $21$ nodes along each side of the computational domain. 
    \label{Fig:StabMDarcy_Five_spot_VMS1_pressure}}
\end{figure}

\begin{figure}
\centering
  \includegraphics[scale=0.45]{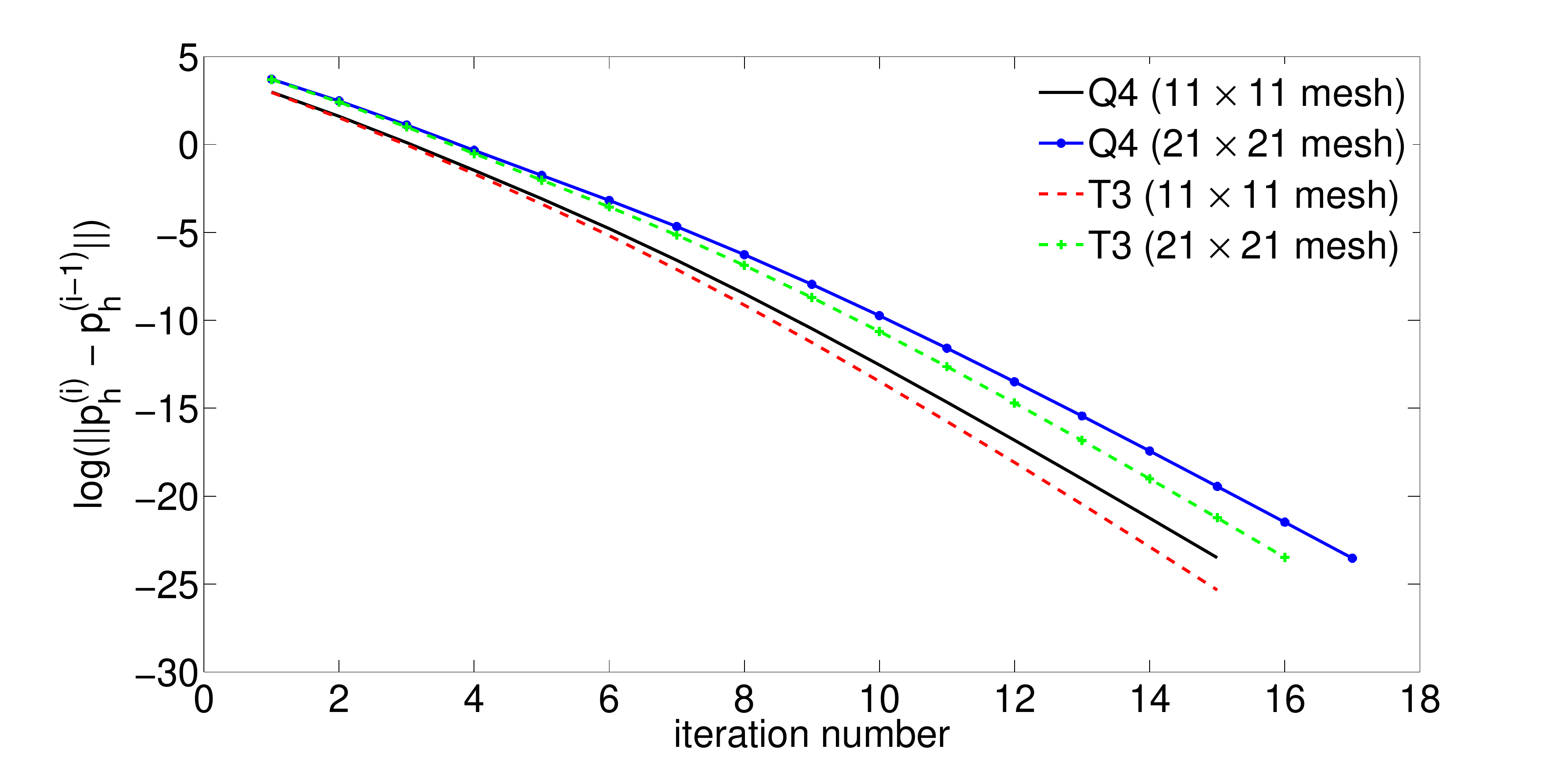}
  \caption{Five-spot problem: Variation of $\|p^{(i)}_h - p^{(i-1)}_h\|$ (which is based on the 
  $2$-norm of the nodal values of the pressure) with respect to iteration number using four-node 
  quadrilateral (denoted by Q4) and three-node triangular (denoted by T3) elements. In this 
  numerical simulation we have used $\epsilon_{\mathrm{TOL}} = 10^{-10}$. 
    \label{Fig:StabMDarcy_Five_spot_VMS1_iteration}}
\end{figure}

\begin{figure}
  \label{Fig:StabMDarcy_five_spot_iters_vs_beta}
  \centering
  \psfrag{beta}{$\beta$}
  \includegraphics[scale=0.37,clip]{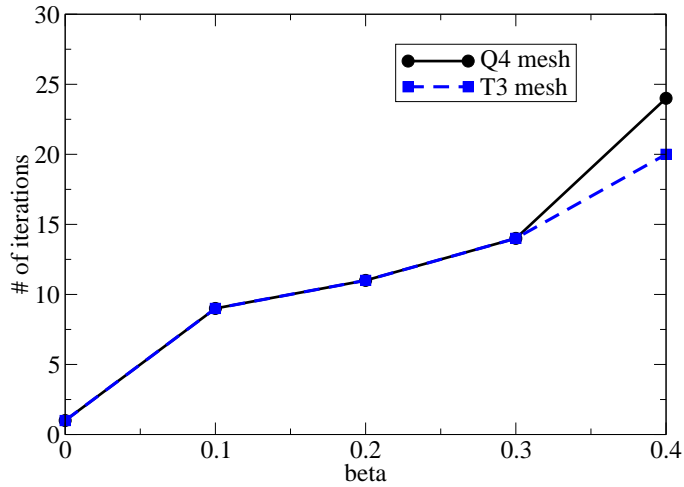}
  \caption{Five-spot problem: This figure shows the number 
    of iterations taken by the numerical algorithm for various 
    values of $\beta$, and for both Q4 and T3 meshes. Note that 
    $\beta$ is an indicator of the strength of the nonlinearity.}
\end{figure}

\begin{figure}
  \label{Fig:StabMDarcy_five_spot_res_vs_iters}
  \centering
  \includegraphics[scale=0.3]{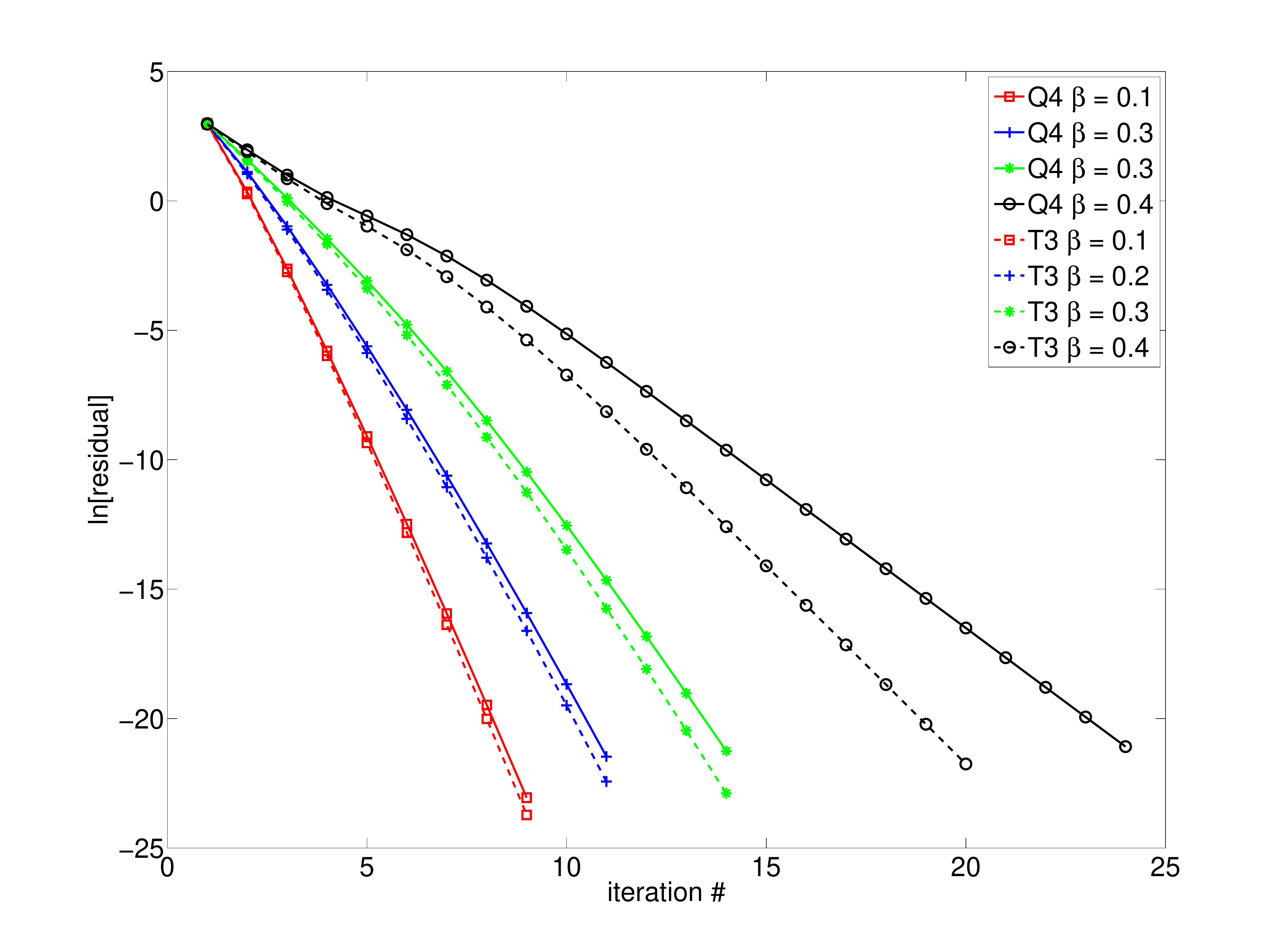}
  \caption{Five-spot problem: This figure shows the 
    variation of the residual with respect to iteration number 
    for various values of $\beta$, and for both Q4 and T3 meshes. 
    Note that the $y$-axis is natural logarithm of the residual.}
\end{figure}

A standard test problem is the quarter five-spot problem. The 
computational domain is a square in the horizontal plane with 
injection and production wells at opposite corners along one 
of the diagonals. The source and sink strengths at injection 
and production wells are, respectively, taken as $+1/4$ and 
$-1/4$. There is no volumetric source/sink (i.e., $\boldsymbol{b}
(\boldsymbol{x}) = \boldsymbol{0}$). The five-spot problem is 
pictorially described in Figure \ref{Fig:StabMDarcy_Five_spot_description}.
Figure \ref{Fig:StabMDarcy_Five_spot_VMS1_pressure} show the 
pressure contours using the proposed stabilized formulation. 
In these simulations we have employed Barus' formula 
\eqref{Eqn:StabMDarcy_non_dimensional_Barus} with $\bar{\alpha_0} 
= 1$, $\bar{\beta} = 0.3$, $\mathcal{A} = 1$, and $\mathcal{B} = 0$. 
The variation of the norm $\|\bar{p}^{(i)}_h - \bar{p}^{(i-1)}_h\|$ 
with respect to iteration number is shown in Figure 
\ref{Fig:StabMDarcy_Five_spot_VMS1_iteration}. As one can see 
from these figures, there are no spurious oscillations in the 
pressure, and the proposed stabilized formulation performed well. 
Figure \ref{Fig:StabMDarcy_five_spot_iters_vs_beta} shows that 
the number of iterations increase with an increase in $\beta$, 
which is a measure of the strength of the nonlinearity. 
Figure \ref{Fig:StabMDarcy_five_spot_res_vs_iters} shows that 
the norm of the residual vector decreases monotonically with 
respect to the iteration number under the proposed stabilized 
formulation for the five-spot problem.

\subsection{The checkerboard problem}
\begin{figure}
\psfrag{L}{$L = 1$}
\psfrag{1}{I}
\psfrag{2}{II}
\psfrag{3}{III}
\psfrag{4}{IV}
\includegraphics[scale=0.75]{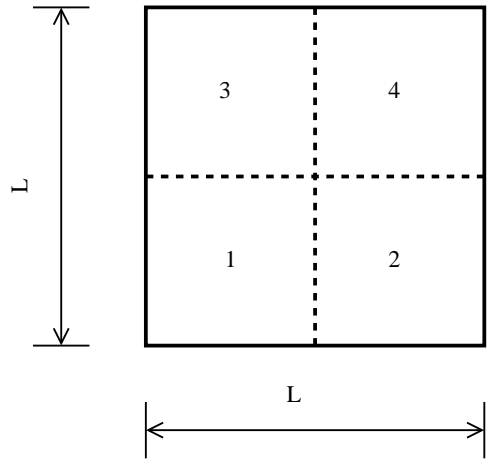}
\caption{Checkerboard problem: A pictorial description. In regions I and IV, we have taken $\bar{\alpha}_0 = 1.0$; and 
in regions II and III, we have taken $\bar{\alpha}_0 = 0.001$. \label{Fig:StabMDarcy_checkered_board_description}}
\end{figure}

\begin{figure}
  \subfigure{
    \includegraphics[scale=0.5]{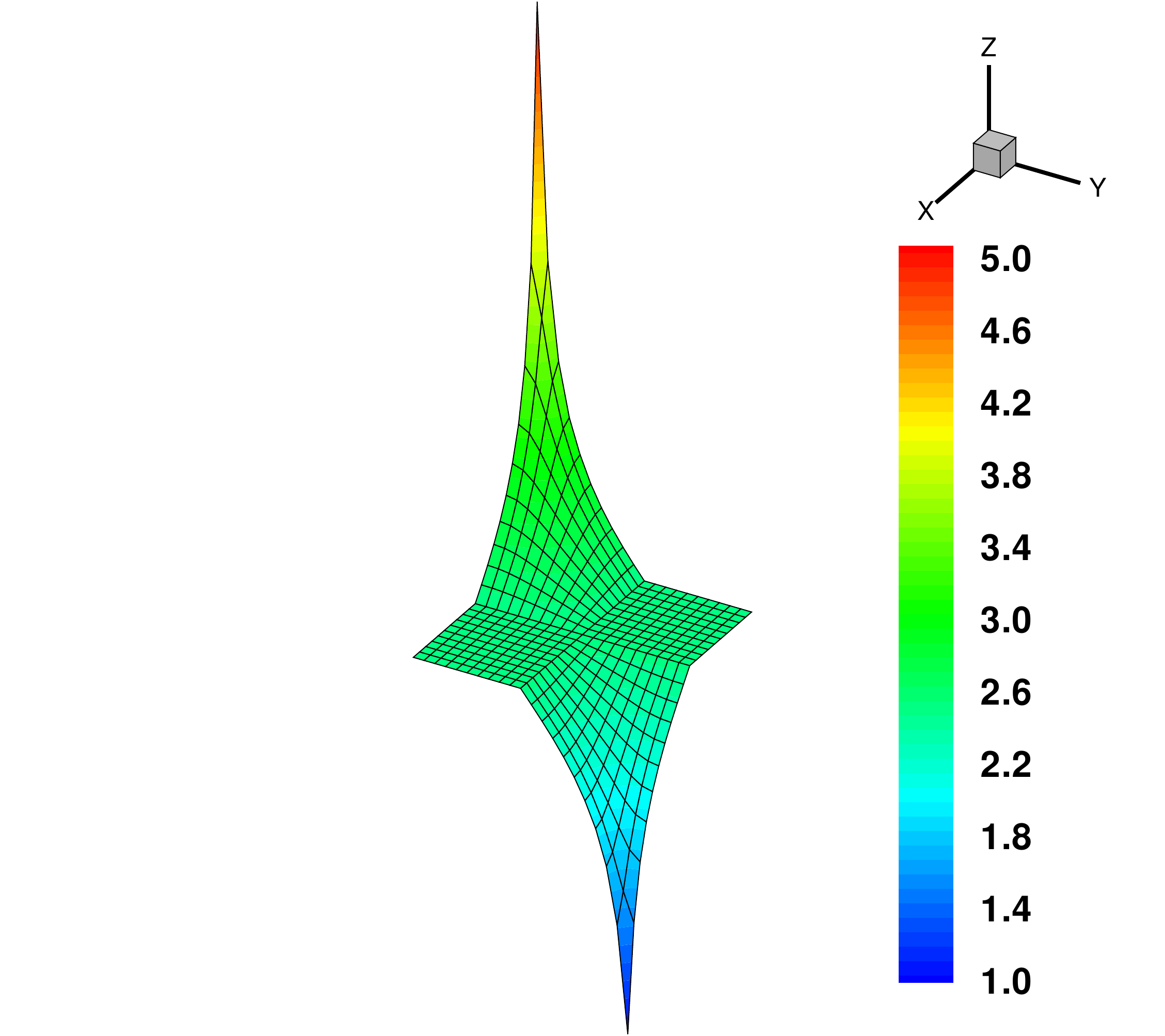}}
  \subfigure{
      \includegraphics[scale=0.5]{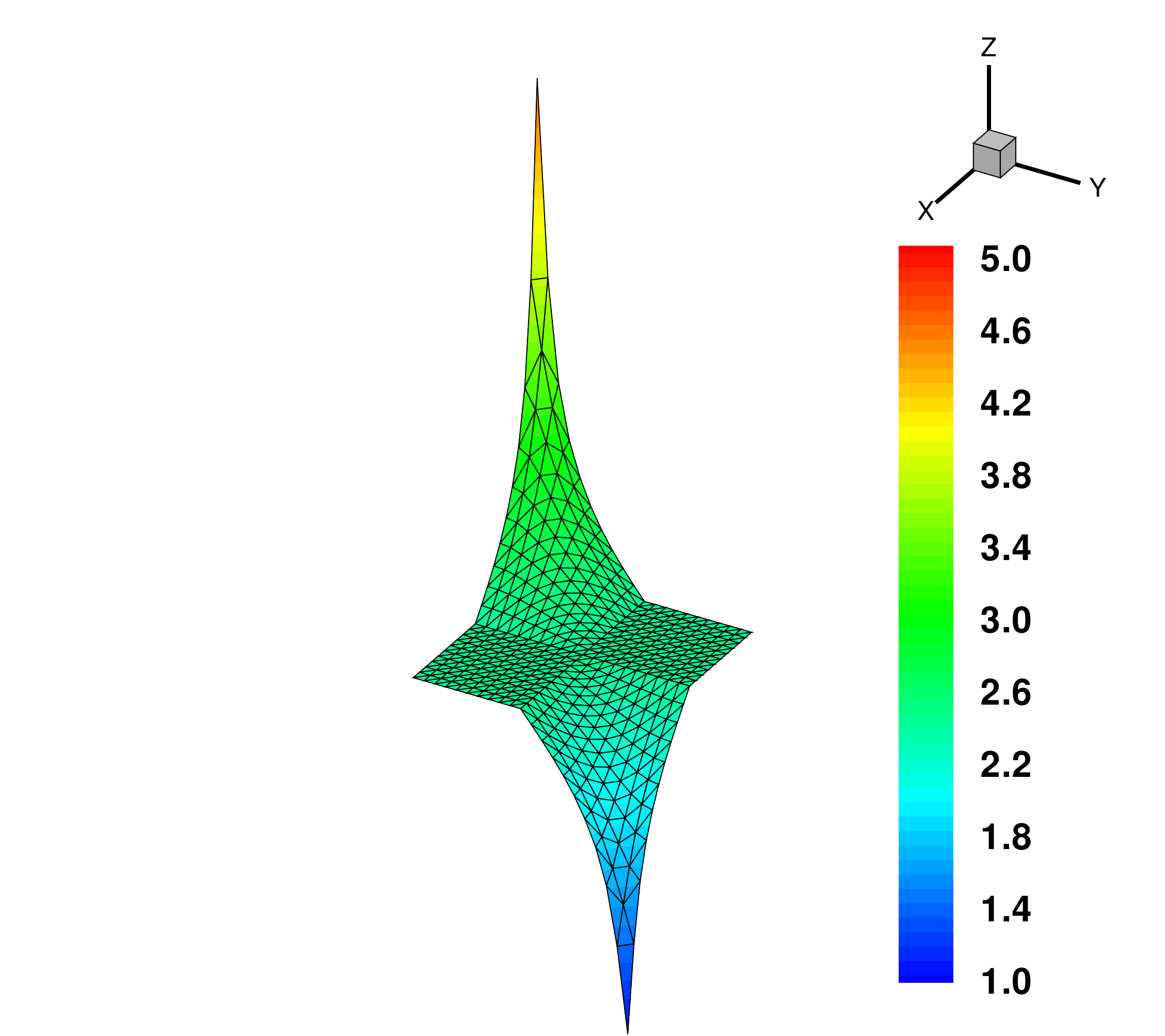}}
  \caption{Checkerboard problem: Pressure contours using four-node quadrilateral (top) and three-node triangular 
  (bottom) elements. We have used $21$ nodes along each side of the computational domain. Barus formula with 
  $\bar{\beta} = 0.3$, $\mathcal{A} = 1$ and $\mathcal{B} = 0$ are employed. \label{Fig:StabMDarcy_Checkerboard_VMS1_pressure}}
\end{figure}

\begin{figure}
\centering
  \includegraphics[scale=0.45]{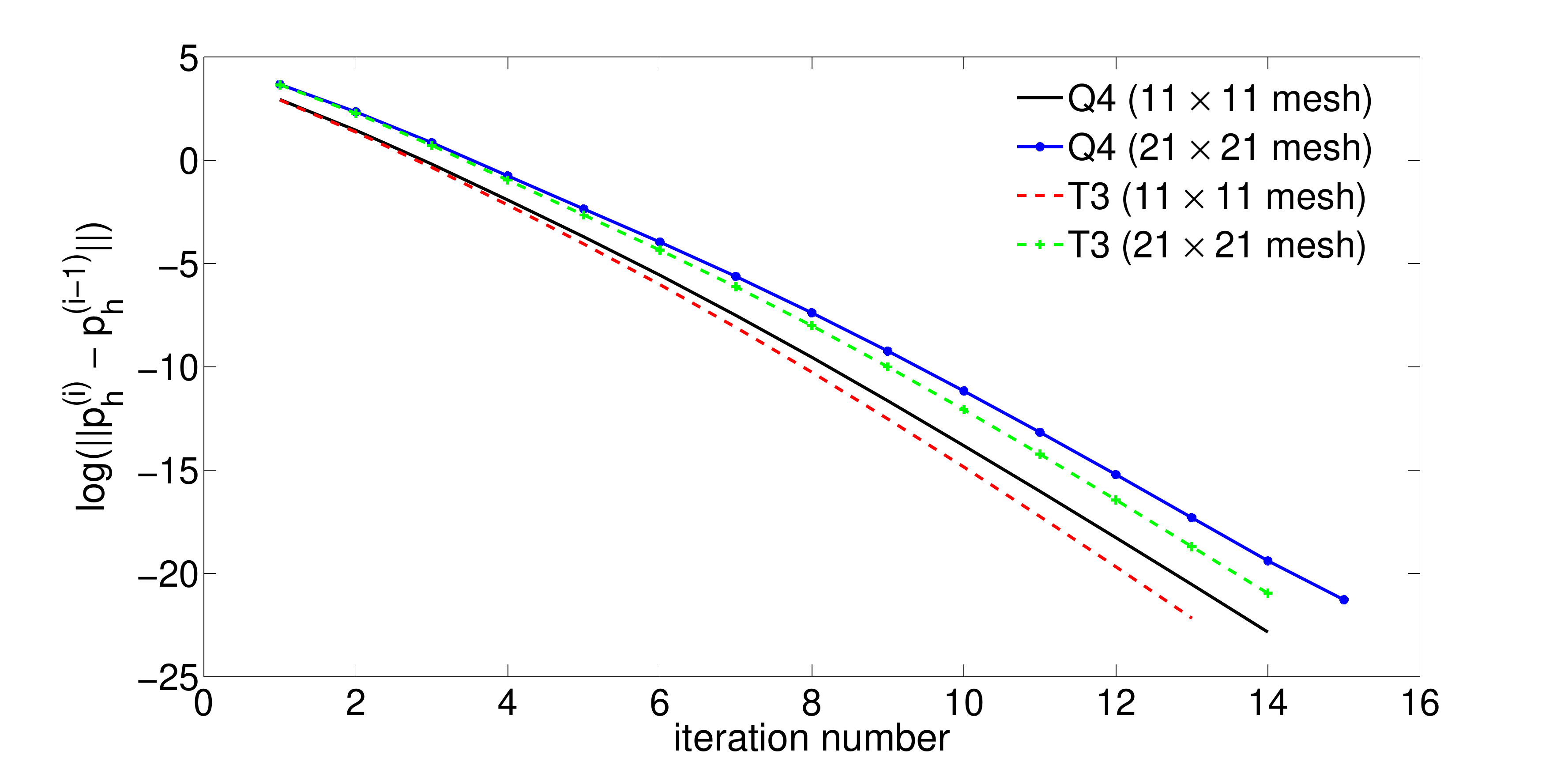}
  \caption{Checkerboard problem: Variation of $\|\bar{p}^{(i)}_h - \bar{p}^{(i-1)}_h\|$ (which is based on the 
  $2$-norm of the nodal values of the pressure) with respect to iteration number using four-node 
  quadrilateral (denoted by Q4) and three-node triangular (denoted by T3) elements. In this 
  numerical simulation we have used $\epsilon_{\mathrm{TOL}} = 10^{-9}$, and $\bar{\beta} = 0.3$. 
    \label{Fig:StabMDarcy_Checkerboard_VMS1_iteration}}
\end{figure}

This problem tests the formulation for the case in which there are 
abrupt changes in the drag coefficient. The geometry and boundary 
conditions are same as the quarter five-spot test problem. But, 
for this test problem, the computational domain is divided into 
four regions as shown in Figure \ref{Fig:StabMDarcy_checkered_board_description}. 
In regions I and IV, we have taken $\bar{\alpha}_0 = 1$; and in regions II 
and III, we have taken $\bar{\alpha}_0 = 0.001$. In the numerical simulations 
we have employed Barus formula with $\bar{\beta} = 0.3$, $\mathcal{A} = 1$, 
$\mathcal{B} = 0$, and have taken $\epsilon_{\mathrm{TOL}} = 10^{-9}$. In Figure 
\ref{Fig:StabMDarcy_Checkerboard_VMS1_pressure} we have shown the contours of 
pressure for both four-node quadrilateral and three-node triangular meshes. In 
Figure \ref{Fig:StabMDarcy_Checkerboard_VMS1_iteration}, we have plotted the 
variation of the norm $\|\bar{p}^{(i)}_h - \bar{p}^{(i-1)}_h\|$ with respect to 
iteration number. As one can see from these figures, there are no spurious 
oscillations in the pressure field, and the proposed stabilized formulation 
performed well. 

\subsection{Numerical $h$-convergence studies}
\begin{figure}
  \includegraphics[scale=0.45]{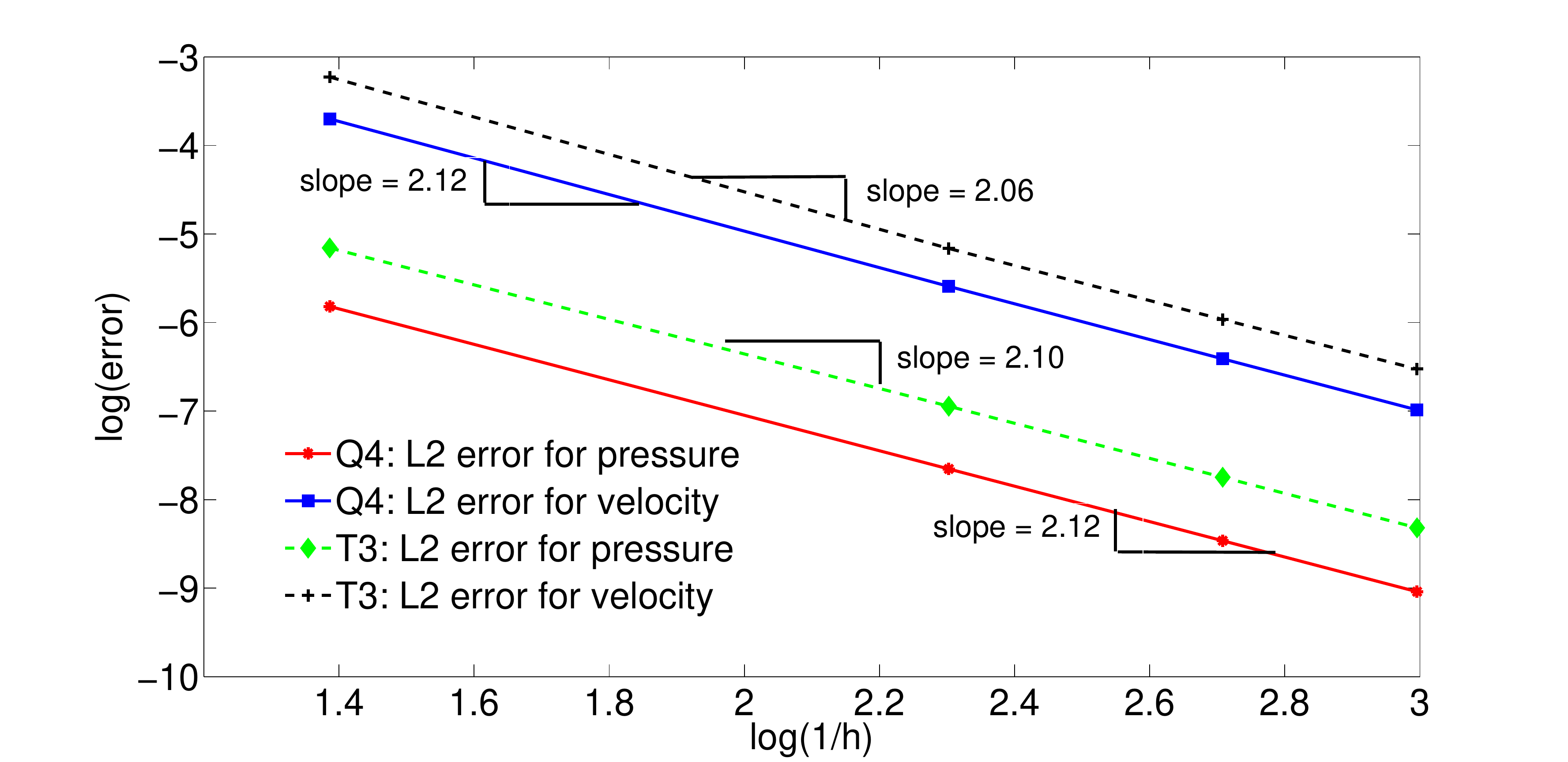}
  \caption{Numerical $h$-convergence studies: The figure shows 
    the rate of convergence in $L_2$-norm for four-node quadrilateral 
    (Q4) and three-node triangular (T3) elements. 
    \label{Fig:StabMDarcy_h_convergence}}
\end{figure}

\begin{figure}
  \includegraphics[scale=0.4]{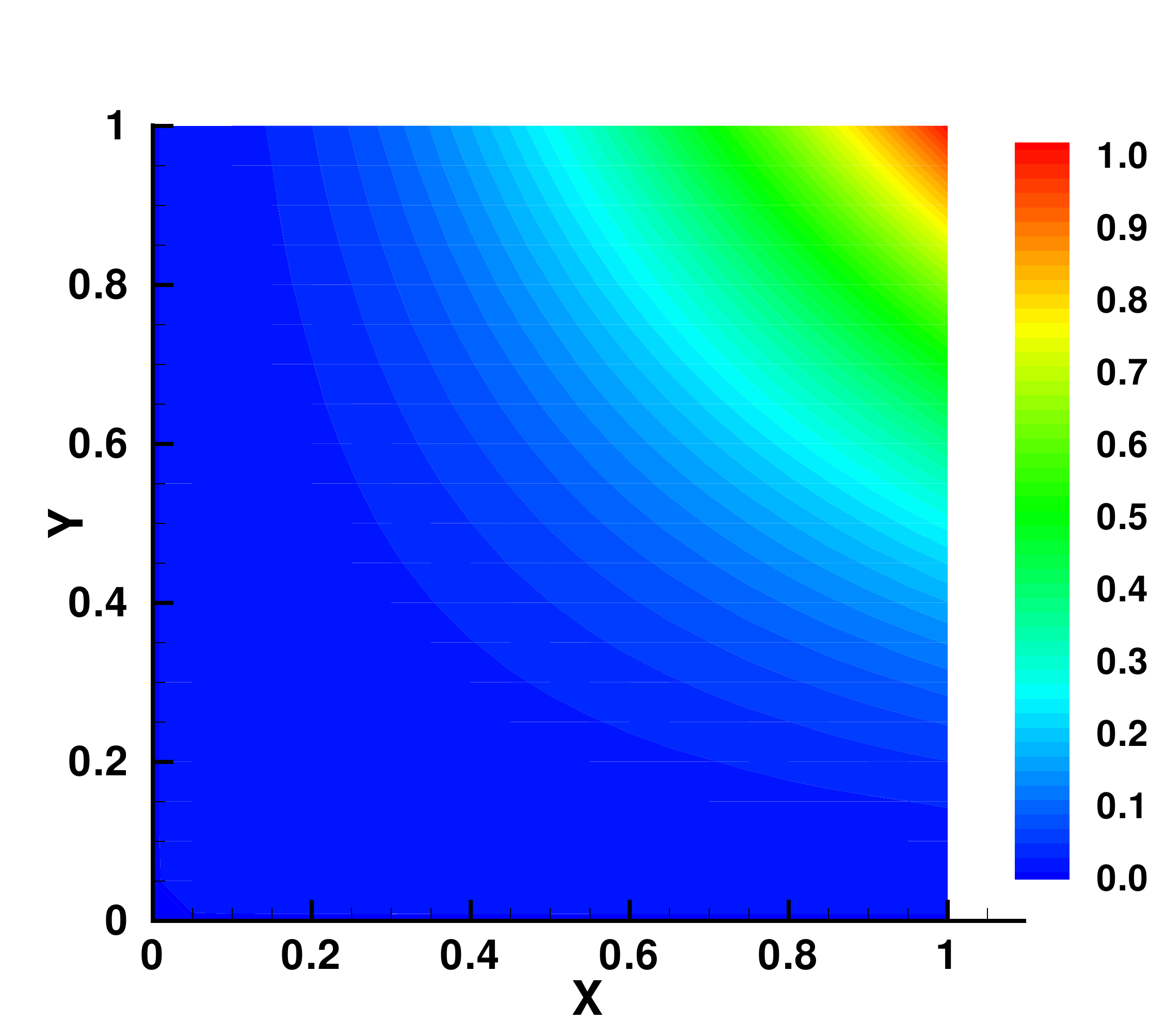}
  \caption{Numerical $h$-convergence studies: Pressure 
    contours using $21 \times 21$ four-node quadrilateral 
    mesh. We have employed equal-order interpolation for 
    the velocity and pressure, and there are no spurious 
    oscillations in the pressure. The analytical solution 
    for pressure is given by $p(x,y) = x^2 y^2$, and the 
    numerical solution matched well with the analytical 
    solution. \label{Fig:StabMDarcy_h_convergence_pressure_contours}}
\end{figure}

Consider the computational domain to be $\Omega = (0,1) \times (0,1)$. We 
have employed Barus formula with $\alpha_0 = 1$ and $\beta = 0.5$. The exact 
solution for the velocity and pressure is given by 
\begin{align}
  v_x(x,y) = +\sin(\pi x) \cos(\pi y), \quad 
  v_y(x,y) = -\cos(\pi x) \sin(\pi y), \quad 
  p(x,y) = x^2 y^2
\end{align}
The density is taken to be unity, and the specific body force is given by 
\begin{align}
  b_x(x,y) = +\exp(x^2 y^2/2) \sin(\pi x) \cos(\pi y) + 2 x y^2, \quad
  b_y(x,y) = -\exp(x^2 y^2/2) \cos(\pi x) \sin(\pi y) + 2 x^2 y
\end{align}
The normal component of the velocity is prescribed to be zero on 
the boundary, and $\epsilon_{\mathrm{TOL}} = 10^{-9}$. The rate of 
convergence with respect to mesh refinement is shown in Figure 
\ref{Fig:StabMDarcy_h_convergence}. The pressure contours are 
shown in Figure \ref{Fig:StabMDarcy_h_convergence_pressure_contours}. 
From these figures it is evident that the proposed algorithm performed 
well.

\subsection{Three-dimensional constant flow}
This problem is a three-dimensional version of the patch test, 
which is employed in Reference 
\cite{Nakshatrala_Turner_Hjelmstad_Masud_CMAME_2006_v195_p4036} 
to assess the stability of a stabilized formulation for Darcy 
equation. Herein, we shall use the patch test to assess the 
stability of proposed stabilized formulation for modified 
Darcy equation. The formulation is considered to pass the 
patch test if the flow field matched with the analytical 
solution up to machine precision.

\begin{figure}
\subfigure{
\includegraphics[scale=0.4]{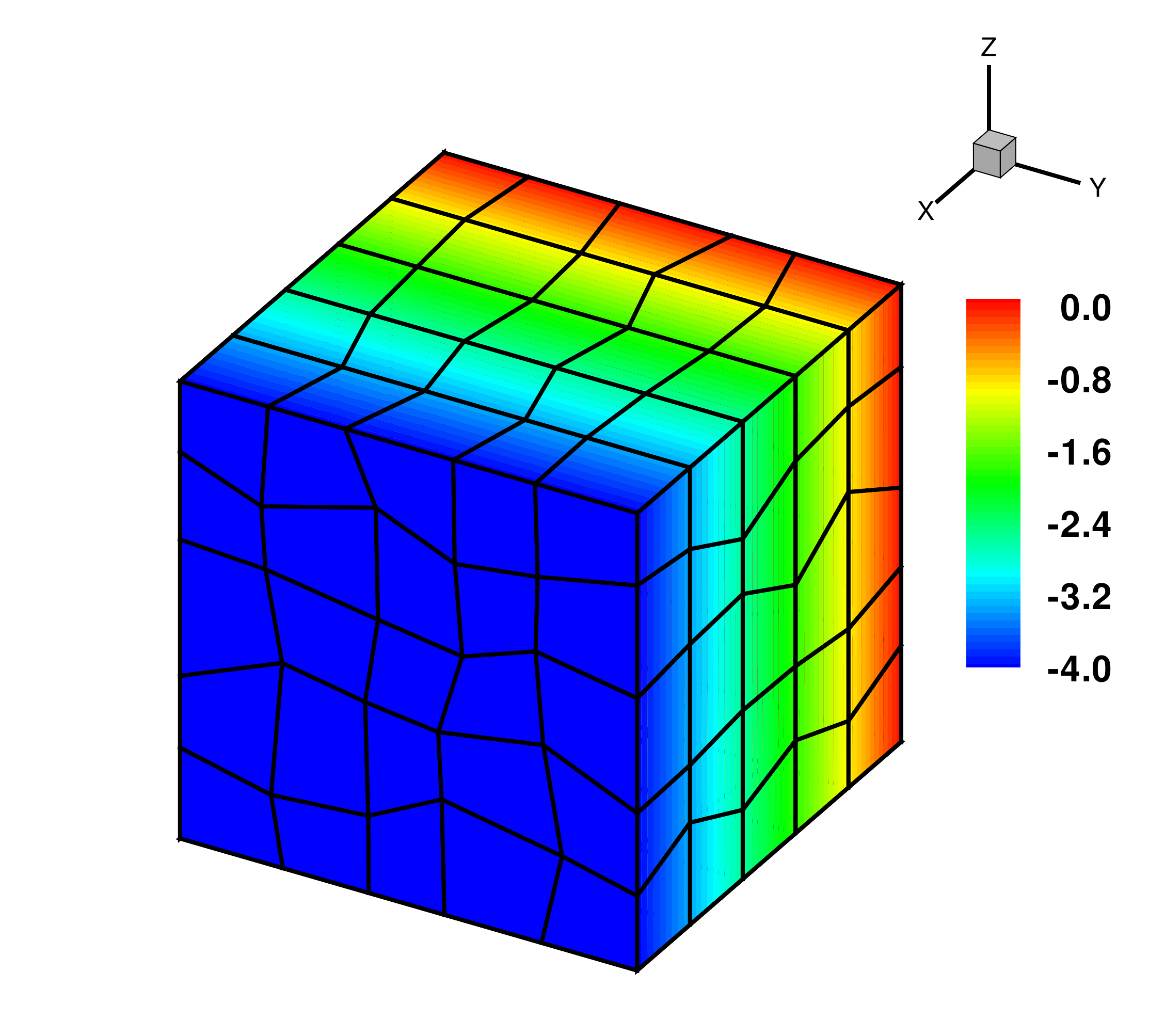}}
\subfigure{
\includegraphics[scale=0.4]{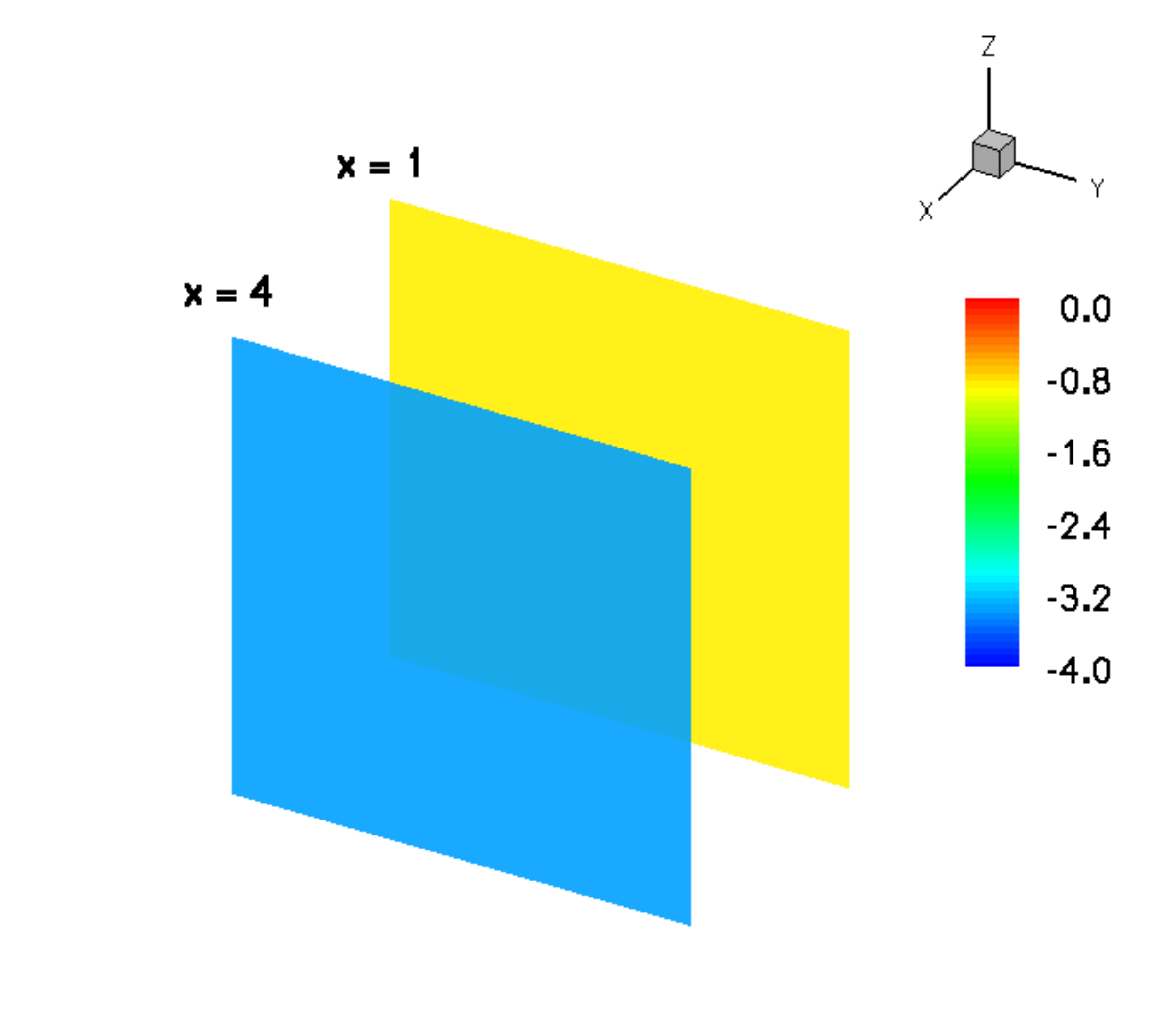}}
\caption{Three-dimensional constant flow: This figure shows the contours of pressure, and the mesh 
is also shown in the top figure. In the bottom figure we have shown the pressure on the $x = 1$ and 
$x = 4$ planes. In this numerical example we have employed Barus formula with $\alpha_0 = 1$ 
and $\beta = 0.1$. \label{Fig:StabMDarcy_3D_patch_test_Q4_pressure}}
\end{figure}

The computational domain is a cube given by $\Omega = (0, 5) 
\times (0, 5) \times (0, 5)$. We have employed Barus formula 
with $\alpha_0 = 1$ and $\beta = 0.1$. We have prescribed the 
normal component of the velocity on $x = 0$ and $x = 5$ faces 
to be unity. On the other four outer faces we have prescribed 
$v_n = 0$. The body force is assumed to be zero, and for 
uniqueness the pressure at the origin is taken to be zero 
(that is, $p_0(0,0,0) = 0$). The analytical solution for 
the pressure is given by
\begin{align}
  p(x,y,z) = -\frac{1}{\beta} \ln[1 + \alpha_0 \beta x]
\end{align}
We have taken $\epsilon_{\mathrm{TOL}} = 10^{-9}$ in this numerical 
simulation. In Figure \ref{Fig:StabMDarcy_3D_patch_test_Q4_pressure} 
we have shown the pressure contours obtained using the proposed 
formulation, and the numerical matched well with the analytical 
solution. 

\subsection{Regions with different permeability}
This problem considers fully saturated, single phase, single 
component flow in the region of a production well close to 
a boundary between regions of permeability that differ by 
orders of magnitude. Figure \ref{fig:Sunset_Geo_mesh} shows 
the problem domain, the production well near the straight 
interface, and the computational mesh employed in the numerical 
simulation. 
The mesh consists of $12,924$ eight-node brick linear elements 
with a total of 63,224 unknowns. The smallest elements (which 
are near the well) are 0.33 units in diameter while the largest 
elements (which are near the perimeter of the circle) are 5.0 
units in diameter. We have taken $\epsilon_{\mathrm{TOL}} = 
10^{-12}$. The properties used in this analysis are listed 
in Table \ref{Table:StabMDarcy_regions_params}. 
The boundary conditions for this problem consists of no-penetration boundaries at the top and bottom of the 
domain (except for the opening at the production well). The pressure at the well opening and around the sides 
is prescribed in a weak fashion with values given in Table \ref{Table:StabMDarcy_regions_params}. The drag 
function $\bar{\alpha}(\bar{p})$ is calculated using equation \eqref{Eqn:StabMDarcy_non_dimensional_linear} 
with $\bar{\beta}$ varying between $0.0$ and $1.0$. 

\begin{figure}[htb!]
\centering
\includegraphics[scale=0.4]{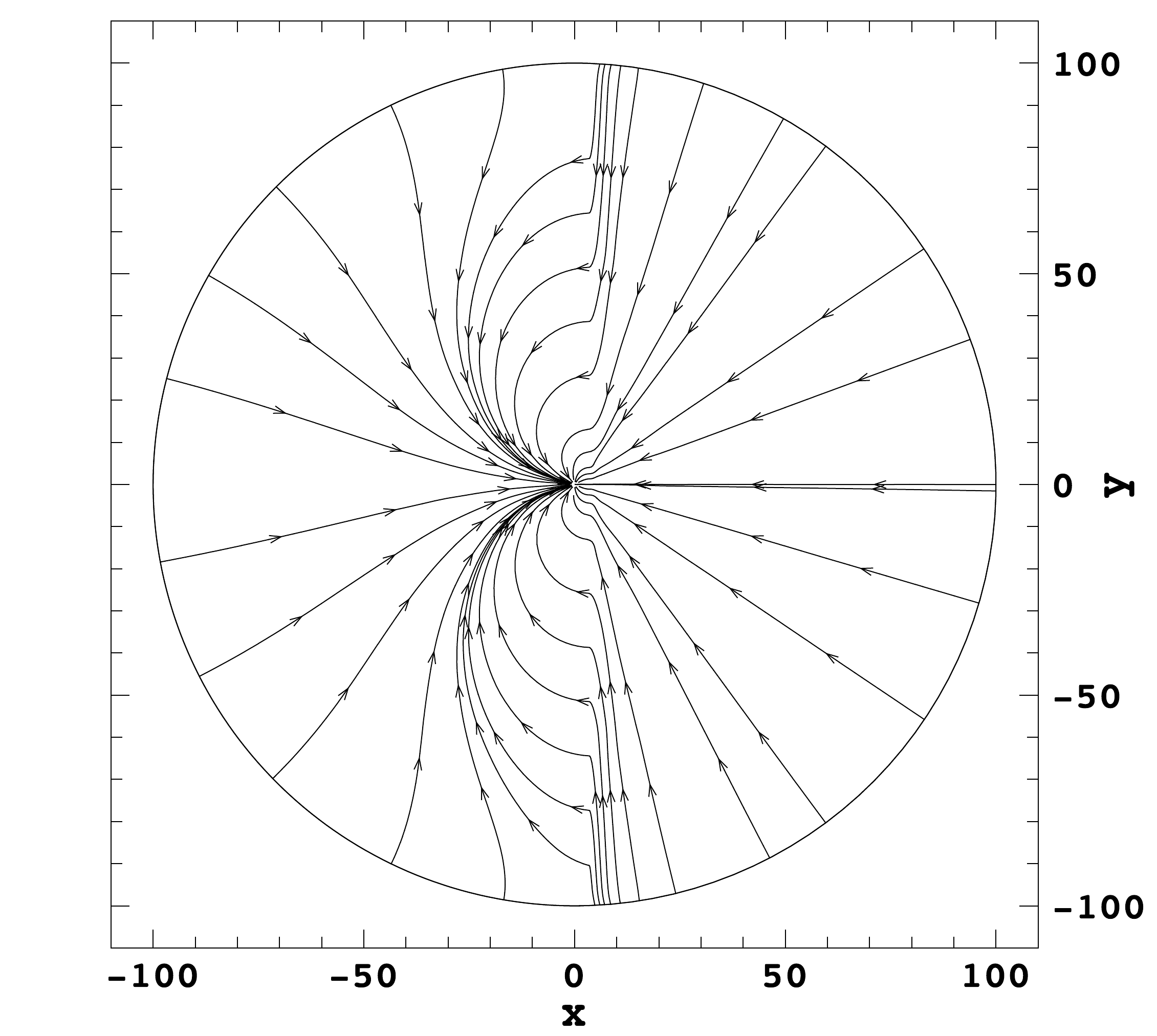}
\caption{Regions with different permeability: This figure shows the velocity streamlines 
for $\beta = 0$. (Only the streamlines for $\beta = 0$ are plotted because the streamlines 
for all other values are similar.)} \label{fig:SunsetStream}
\end{figure}

\begin{figure}[htb!]
  \centering
  \includegraphics[scale=0.3]{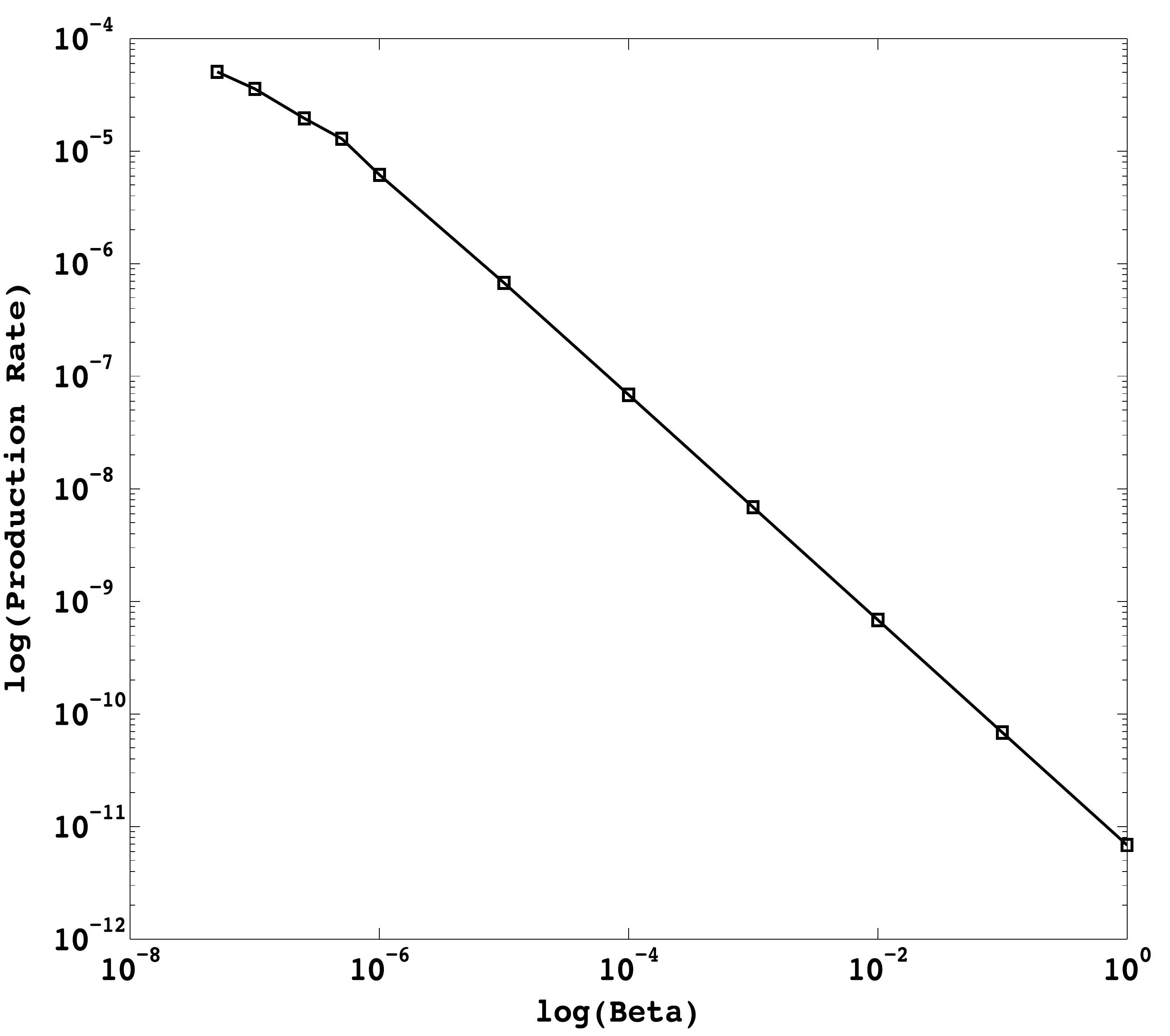}
  \caption{Regions with different permeability: This figure shows 
    the production rate at the well opening for various values of 
    $\beta$.} \label{fig:SunsetProd}
\end{figure}

\begin{figure}[htb!]
  \centering
  \includegraphics[scale=0.3]{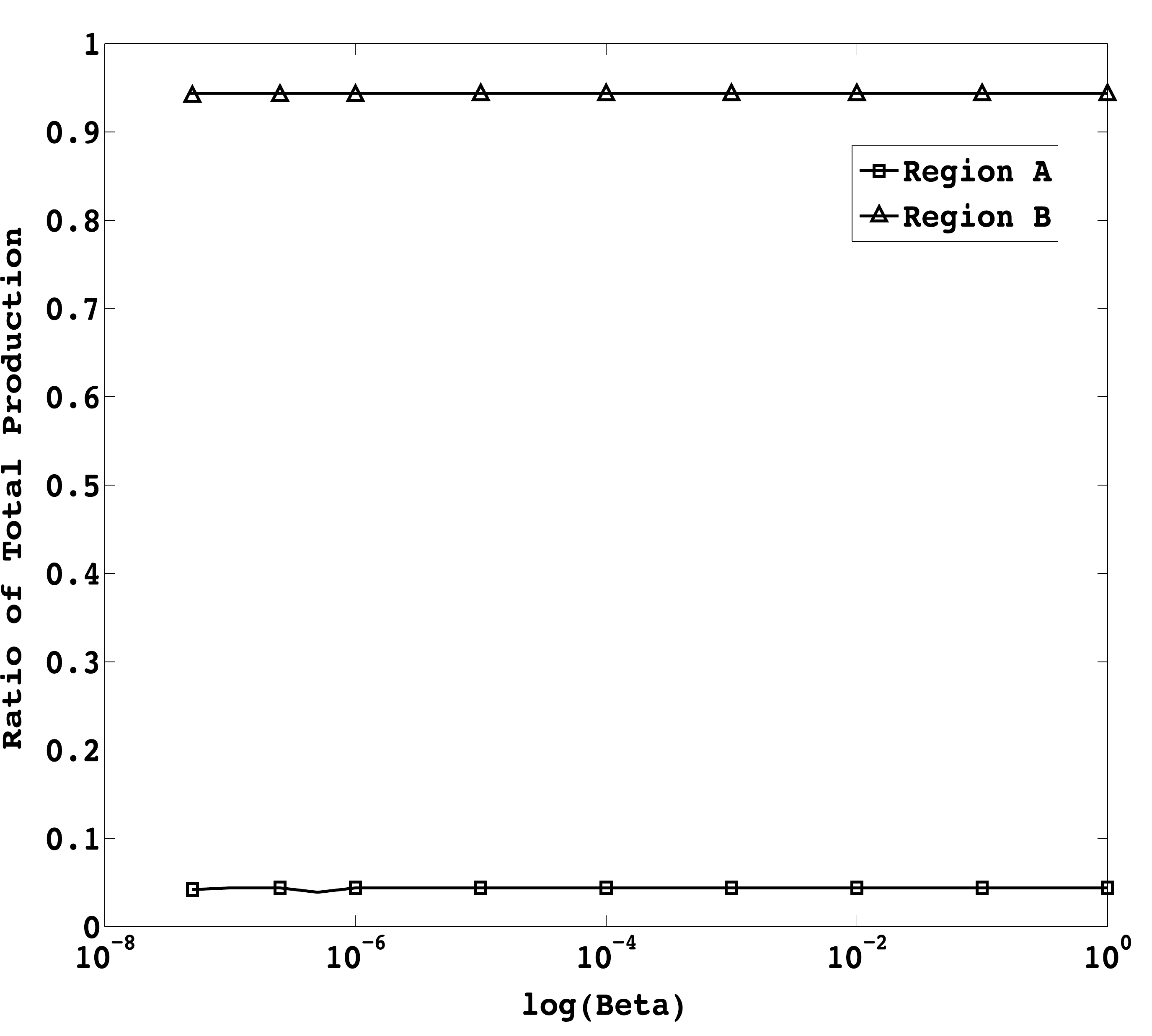}
  \caption{Regions with different permeability: This figure shows the 
    ratio of the total production emanating from regions A and B with 
    respect to $\beta$.} \label{fig:SunsetRatio}
\end{figure}

The exact solution for the velocity streamlines for this problem 
for constant drag function can be obtained by the method of images 
and is presented in Reference \cite{Bear_Hydraulics}. In region A, 
the streamlines should curve away from the  production well in a 
radial fashion. In region B, the streamlines should remain straight. 
Figure \ref{fig:SunsetStream} shows similar results for the streamlines 
as produced by the method presented in this work. Note that the streamlines 
in Region B are straight as predicted by the exact solution. Figure 
\ref{fig:SunsetProd} shows a log-log plot of the production rate at 
the well as a function of $\bar{\beta}$. The production rate was 
computed as the total flux of fluid $\int \boldsymbol{v}(\boldsymbol{x}) 
\cdot \boldsymbol{n}(\boldsymbol{x}) \ \mathrm{d}\Gamma$ across the well 
opening. 
\emph{Notice that as $\bar{\beta}$ increases, the production rate decreases in a linear fashion. Figure \ref{fig:SunsetRatio} 
shows the amount of the total production at the well that emanates from region A or B. This result is of particular interest in 
determining the effect of $\bar{\beta}$ on the amount which the well will draw from either region. The results show that 
although as $\bar{\beta}$ increases the total production decreases, the proportion of the total production that comes 
from region A or B remains constant.}

\begin{table}[ht!]
\caption{Parameters used in regions with different permeability problem}
\centering
\begin{tabular}{ll}
\hline
Parameter & Value  \\
\hline
$\bar{\alpha}_0$ region A  & 0.001  \\
$\bar{\alpha}_0$ region B  & 1.0      \\
$\bar{\rho}$ & 1.0 \\
$\bar{\boldsymbol{b}}$ & (0,0,0) \\
$\bar{p}_0$  far-field boundary & 1.0 \\
$\bar{p}_0$  well opening & 3.333 $\times 10^{-4}$  \\
Radius of far-field boundary & 100.0 \\
Radius of well & 1.0 \\
\hline
\end{tabular}
\label{Table:StabMDarcy_regions_params}
\end{table}

\subsection{$\mathrm{CO}_2$ leakage through an abandoned well (large-scale problem)}
\begin{figure}[htb!]
  \centering
  \includegraphics[scale=0.4]{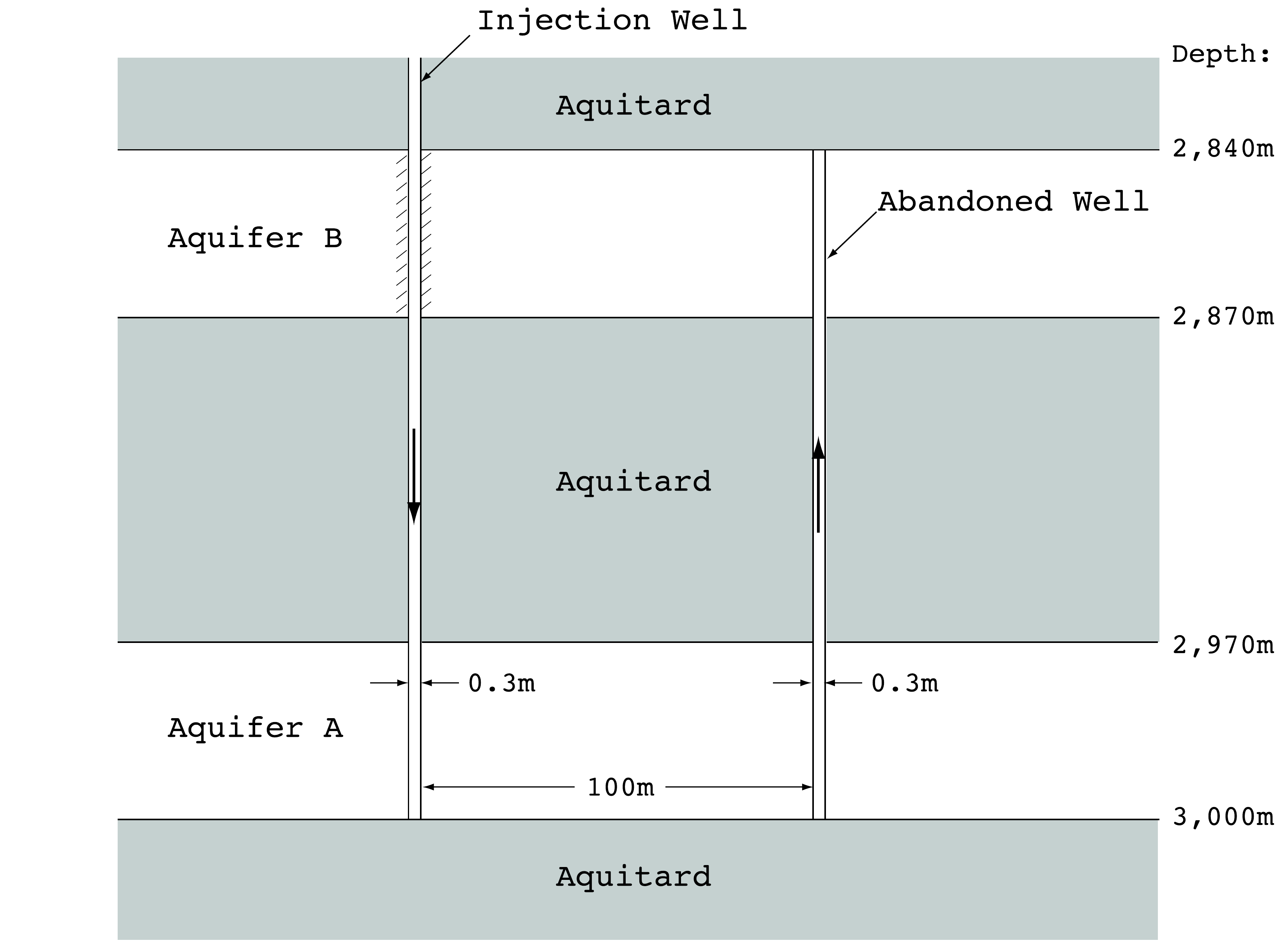}
  \caption{$\mathrm{CO}_2$ leakage through an abandoned well: A pictorial 
    description of the problem showing the cross-sectional view along the plane 
    containing the injection and abandoned wells.} \label{fig:LeakyWellDomain}
\end{figure}

\begin{figure}[htb!]
  \centering
  \subfigure{\includegraphics[scale=0.3]{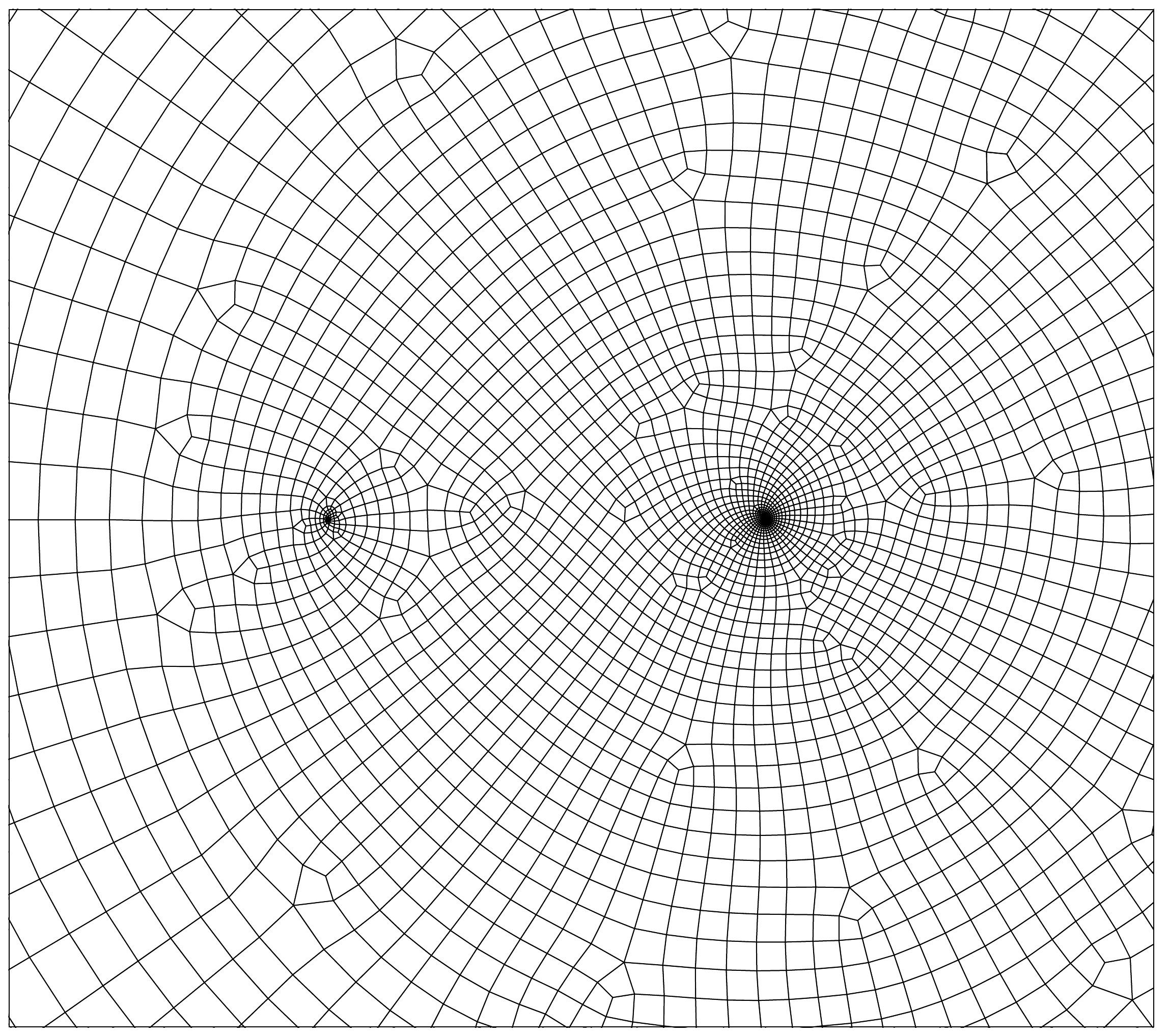}}
  \subfigure{\includegraphics[scale=0.3]{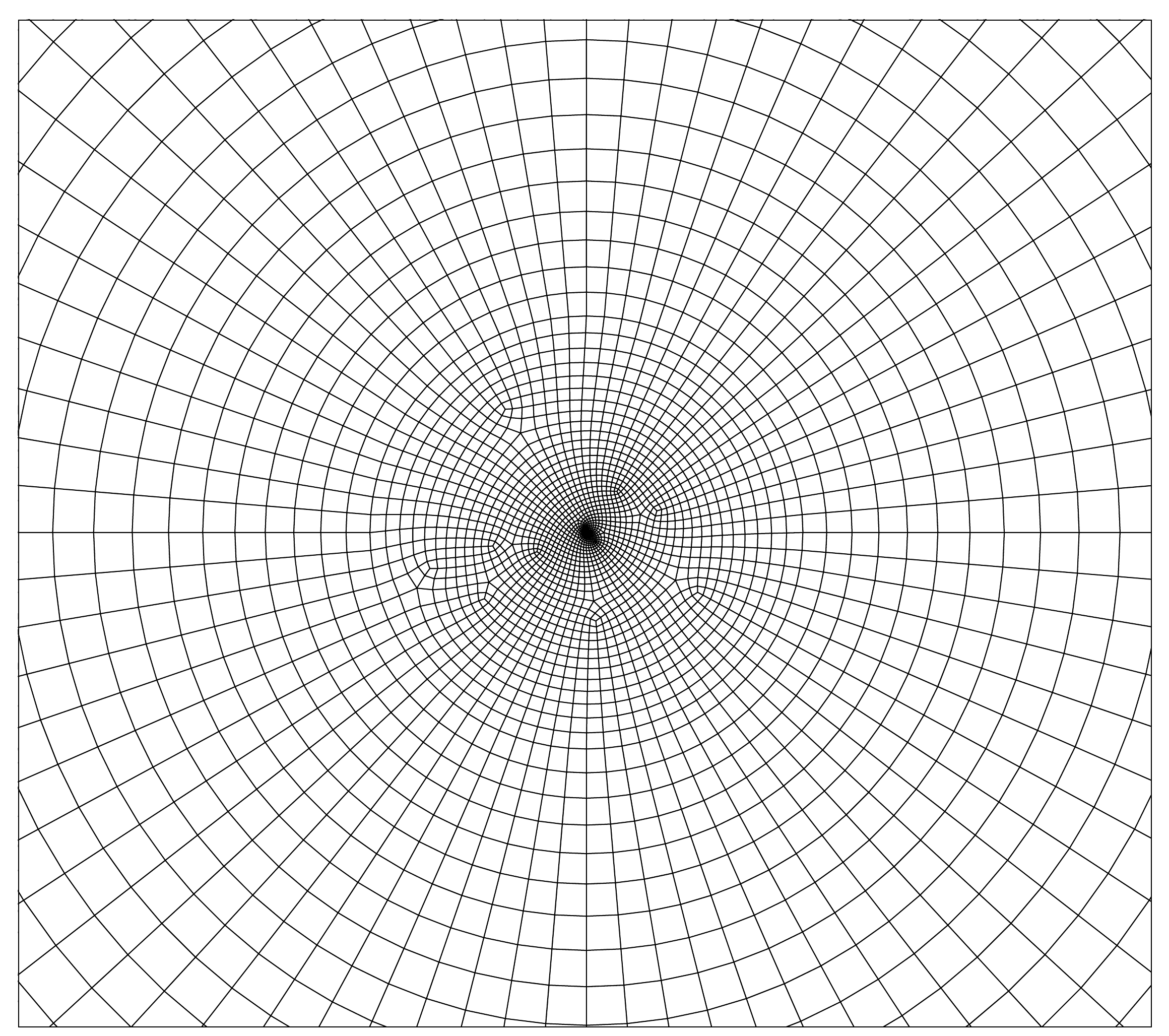}}
  \subfigure{\includegraphics[scale=0.7]{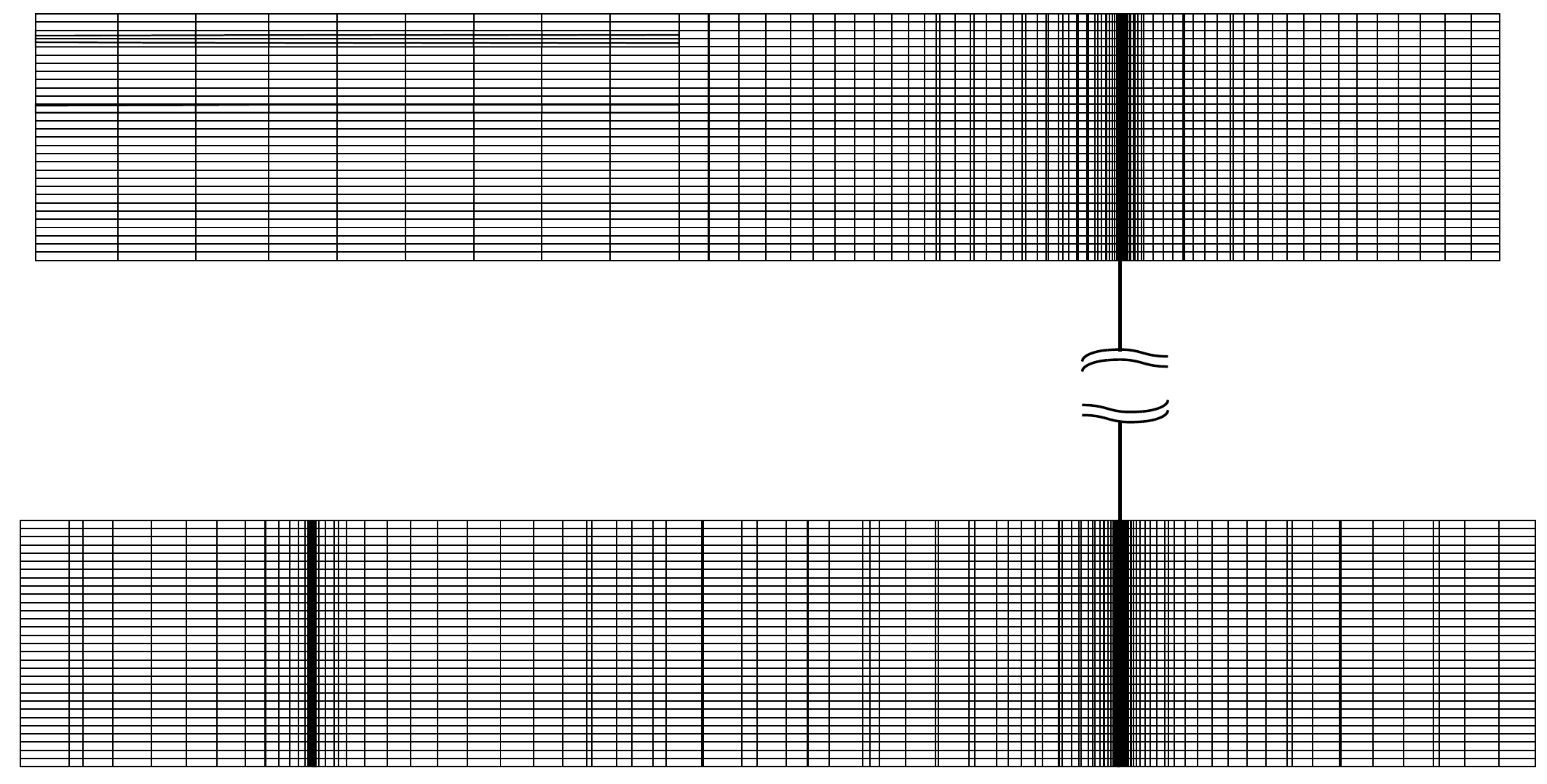}}
  \caption{$\mathrm{CO}_2$ leakage through an abandoned well: This figure 
  shows the computational mesh used in the numerical simulation. The top-left 
  subfigure shows the mesh at the bottom of the domain, and the  top-right 
  subfigure shows the mesh at the top of the domain. The bottom subfigure 
  shows the slice through the center near the injection and abandoned wells. 
  There are over $1.14$ million unknowns in this test problem.}\label{fig:LeakComputationalMesh}
\end{figure}

The last numerical example makes an important contribution to the 
study of geological carbon-dioxide ($\mathrm{CO}_2$) sequestration 
into underground aquifers. When $\mathrm{CO}_2$ is pumped into an 
underground aquifer, it can leak through fissures in the surrounding 
aquitard layers or though man-made penetrations, such as abandoned 
wells. An important question regarding the suitability of certain 
locations for sequestration involves predicting the leakage rate of 
$\mathrm{CO}_2$ into other aquifers. This  numerical example models 
the leakage of $\mathrm{CO}_2$ through an abandoned well as it 
is injected into an aquifer. The geometry is shown in Figure 
\ref{fig:LeakyWellDomain}. The diameter of the wells are $0.3 
\; \mathrm{m}$ and they are located $100 \; \mathrm{m}$ apart 
from each other. The far-field boundary is located at a radius 
of $500 \; \mathrm{m}$. The computational mesh for this problem, 
which shown in Figure \ref{fig:LeakComputationalMesh}, consists 
of 271,050 eight-node brick linear elements and 284,019 nodes 
for a total of 1.14 million unknowns. Near the injection well 
and the abandoned well the element diameter is $0.05 \; \mathrm{m}$. 
Near the far-field boundary the element diameter is $2.5 \; \mathrm{m}$. 
To simplify the problem we assume fully saturated, single component, 
single phase incompressible flow and that there is no body force, 
$\bar{\boldsymbol{b}} = (0,0,0)$. The boundary conditions consist 
of no penetration boundaries on the top and bottom of the aquifers 
and weakly prescribed pressure boundary conditions along the sides. 
For aquifer B, a constant reference pressure, $p_{\mathrm{ref}}$, of 
$2.9315 \times 10^7$ Pa is prescribed along the outer boundary. For 
aquifer A, a constant reference pressure of $3.0599 \times 10^7$ Pa 
is prescribed. At the injection well, a constant inflow velocity of 
$0.262$ m/s is prescribed. The parameters used for this problem are 
listed in Table \ref{table:Probelm2Params}. $\bar{\alpha}(\bar{p})$ 
is determined according to the exponential function given in equation 
\eqref{Eqn:StabMDarcy_non_dimensional_Barus}, with $\bar{\beta}$ varying 
from 0.0 to 1.0. We have again taken $\epsilon_{\mathrm{TOL}} = 10^{-12}$. 

This problem was solved in parallel setting with 32 processors, on the 
Tri-Lab Linux Capacity Cluster, using Aria computer code \cite{Aria}. 
Load balancing across the processors was achieved using the Zoltan 
\cite{Devine_Boman_Heaphy_Hendrickson_Vaughan_ComputSciEng_2002_v4_p90} 
package. The total CPU time for a single run with a given $\beta$ was 
35 minutes with an average memory use on each processor of 400 MB. A 
complete investigation into the parallel scalability of this algorithm 
is intended for future work.

\begin{figure}[htb!]
  \centering
  \subfigure{\includegraphics[scale=0.3]{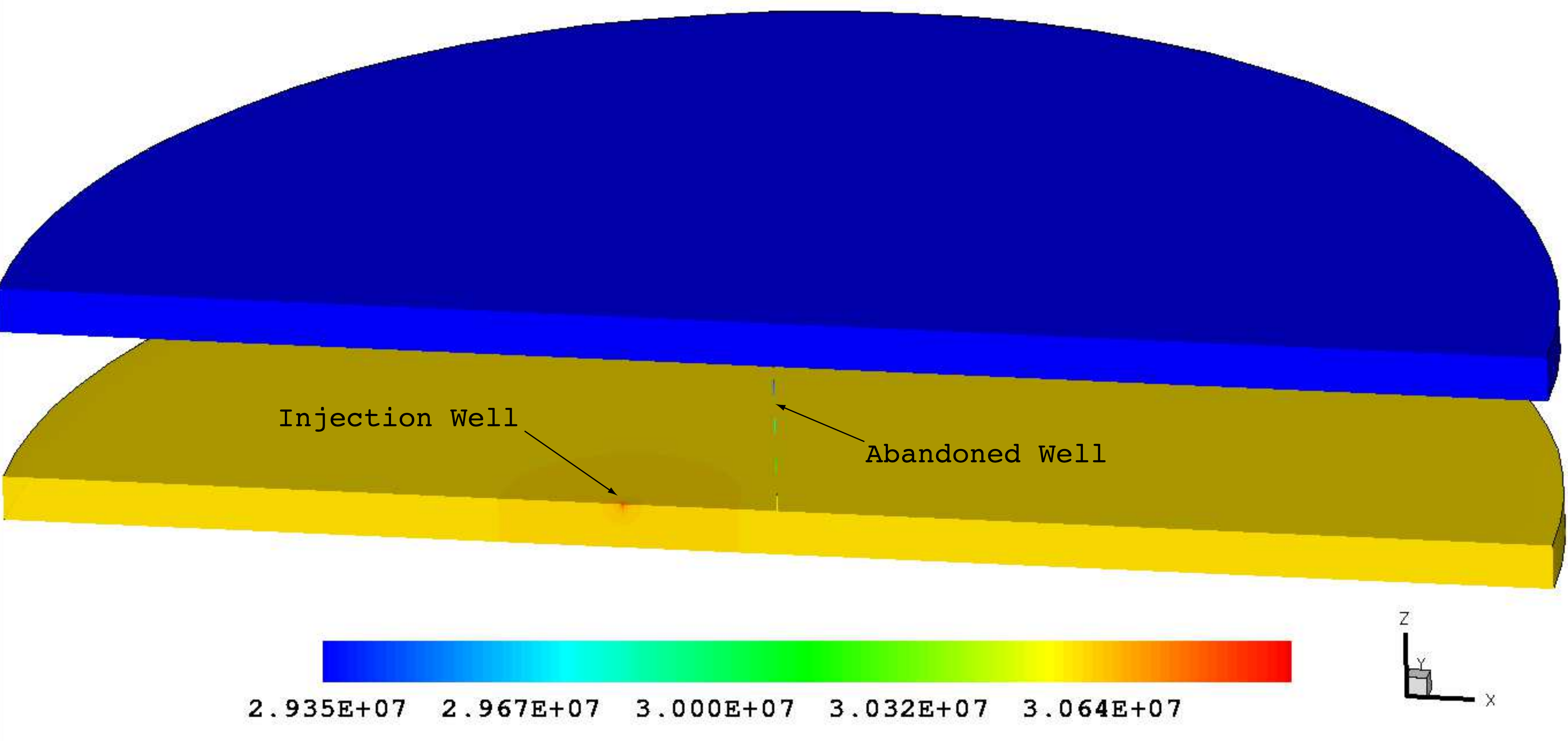}}
  \subfigure{\includegraphics[scale=0.3]{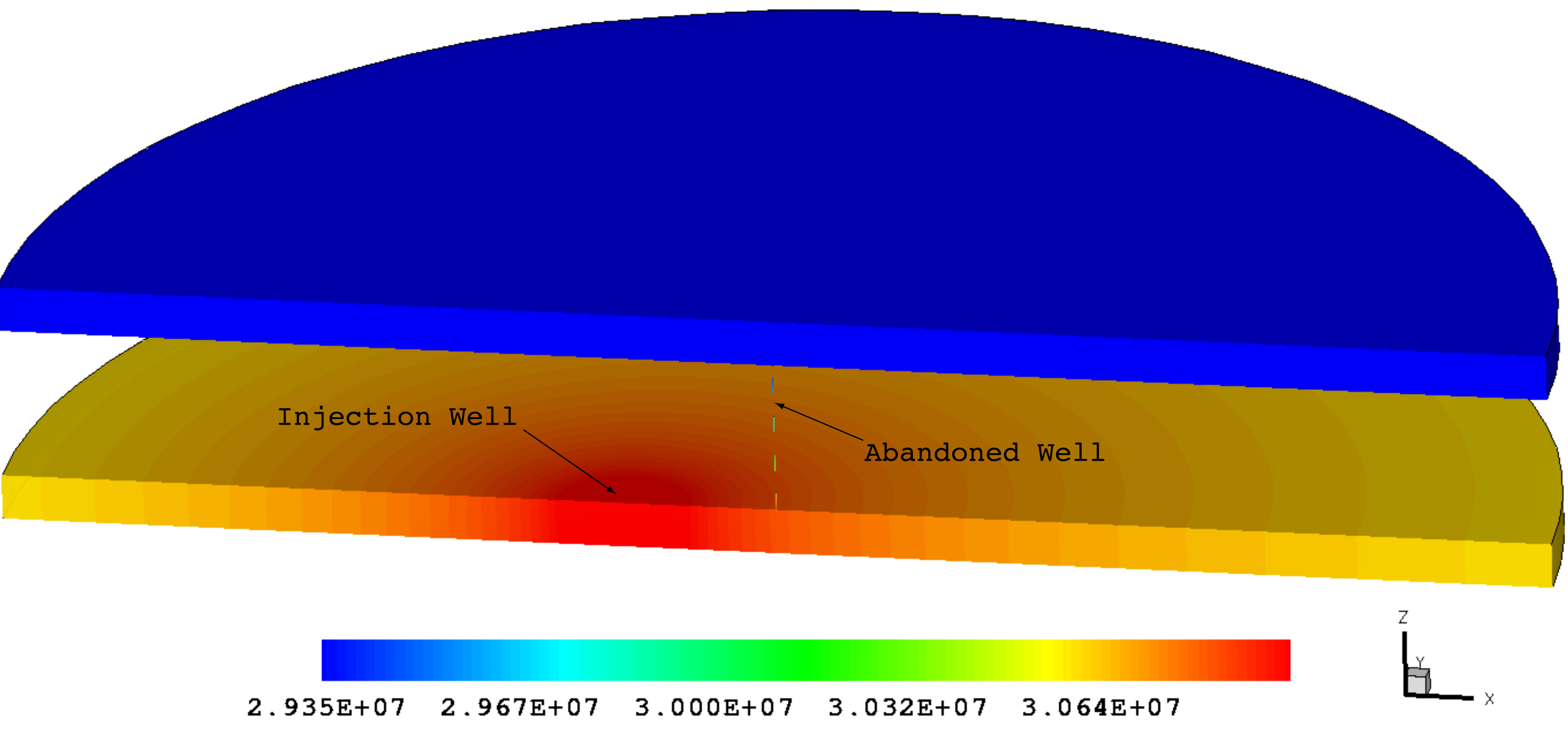}}
  \caption{$\mathrm{CO}_2$ leakage through an abandoned well: 
    This figure shows contours of the pressure for $\bar{\beta} 
    = 0$ (top) and $\bar{\beta} = 1$ (bottom). It is evident from 
    the figure that the modified Darcy model predicts higher 
    pressures and higher pressure gradients than of Darcy model. 
    Half of the domain has been removed to show the detail.} 
  \label{fig:LeakPressure}
\end{figure}

\begin{figure}[htb!]
  \centering
  \subfigure{\includegraphics[scale=0.3]{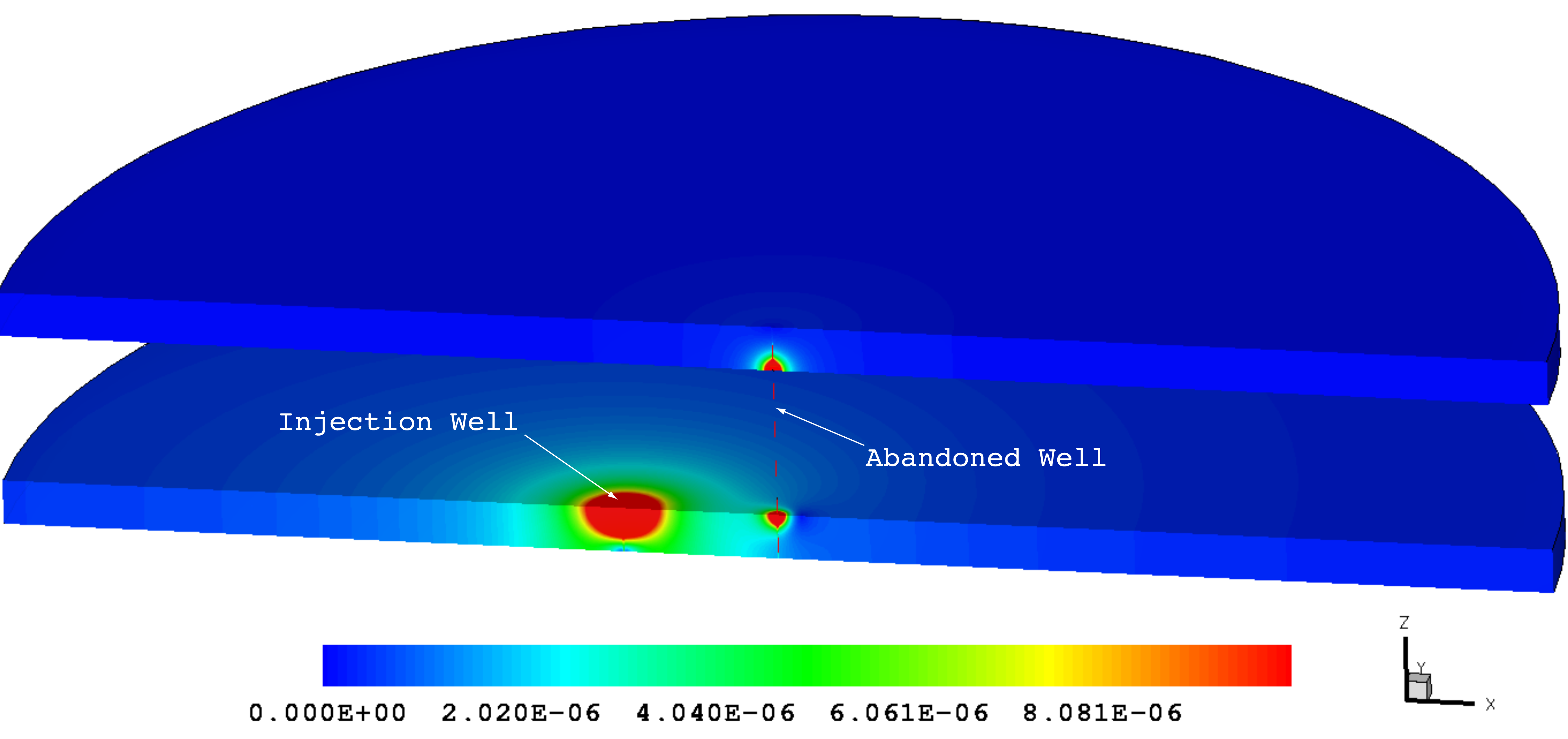}}
  \subfigure{\includegraphics[scale=0.3]{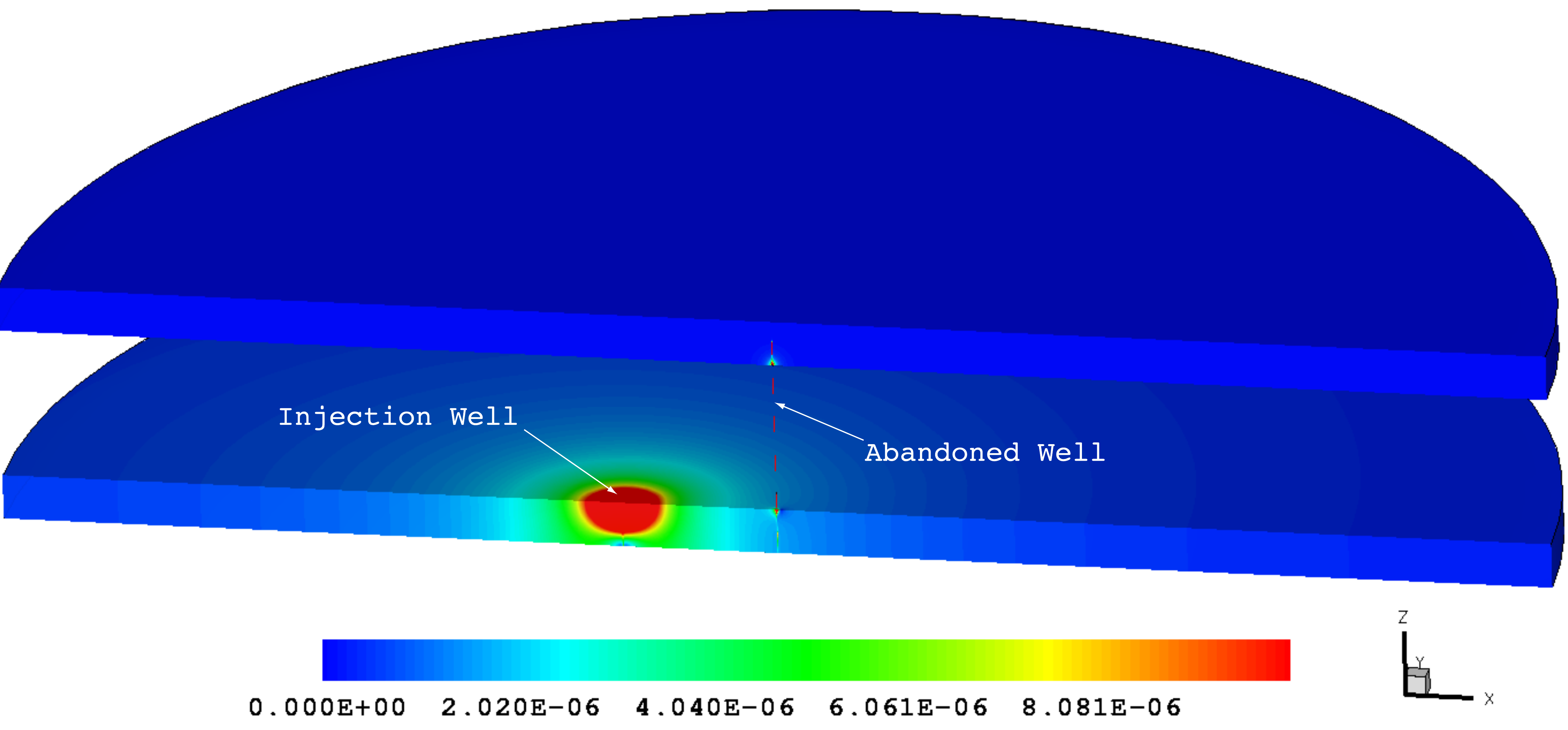}}
  \caption{$\mathrm{CO}_2$ leakage through an abandoned well: This 
    figure shows contours of the magnitude of the velocity for 
    $\bar{\beta} = 0$ (top) and $\bar{\beta} = 1$ (bottom). 
    It is evident from the figure that Darcy model predicts 
    higher velocity in the abandoned well than the prediction 
    made by the modified Darcy model. Half of the domain has 
    been removed to show the detail.}\label{fig:LeakVel}
\end{figure}

\begin{figure}[htb!]
  \centering
  \includegraphics[scale=0.3]{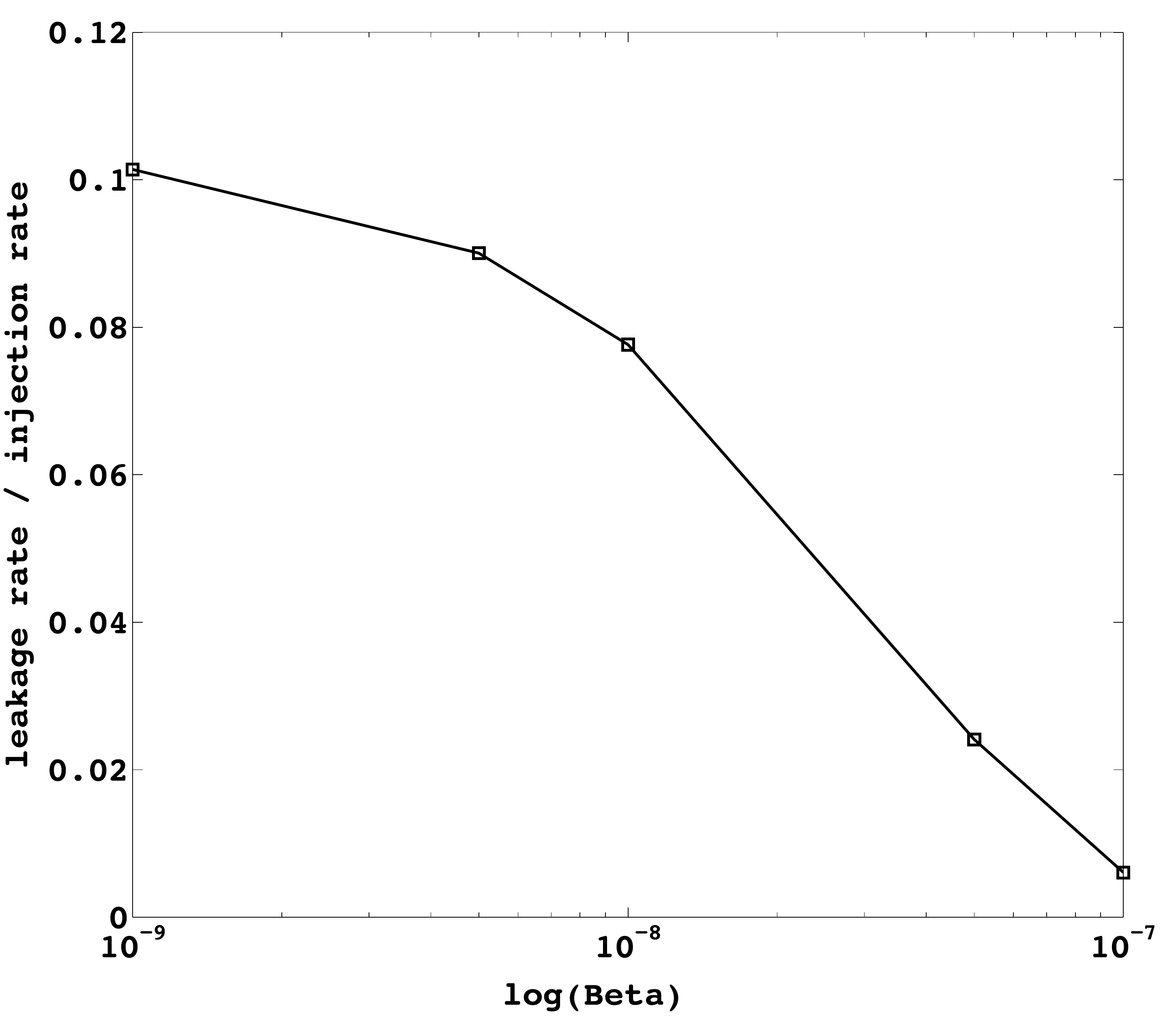}
  \caption{$\mathrm{CO}_2$ leakage through an abandoned well: 
    This figure shows the ratio of leakage rate to injection 
    rate versus $\beta$. Darcy model over-predicts the leakage 
    rate than the modified Darcy model, and the ratio between 
    the leakage rate and the injection rate decreases monotonically 
    with increase in $\beta$. This discovery will be crucial in 
    designing the seal in a geological carbon-dioxide sequestration 
    geosystem. A design of the seal based on Darcy model will be 
    too conservation, which may increase the expense of the 
    geological carbon-dioxide sequestration project.} 
  \label{fig:LeakageInjection}
\end{figure}

Figure \ref{fig:LeakPressure} shows contours of the pressure for the minimum and maximum values of $\bar{\beta}$. Notice 
that as $\bar{\beta}$ increases, the pressure in the region surrounding the injection well increases. Figure \ref{fig:LeakVel} 
shows the magnitude of the velocity for the maximum and minimum values of $\bar{\beta}$. As $\bar{\beta}$ increases, the 
velocity magnitude in the region of the abandoned well decreases suggesting that not including the pressure dependent viscosity 
under predicts the leakage rate. Figure \ref{fig:LeakageInjection} shows the ratio of injection rate to leakage rate as $\beta$ 
increases. The injection rate was computed as the total flux of fluid $\int \boldsymbol{v}(\boldsymbol{x}) \cdot \boldsymbol{n}(\boldsymbol{x}) \ \mathrm{d}\Gamma$ 
across the injection well opening. Likewise, the leakage rate was computed as the total flux of fluid across a cross-section of 
the abandoned well at the bottom of the aquifer B.  In References \cite{Ebigbo_Class_Helmig_ComputGeoscience_2007_v11_p103,
Nordbotten_Celia_Bachu_Dahle_EnvSciTech_2005_v39_p602}, the authors present transient results for a 
similar simulation for multiphase flow of $\mathrm{CO}_2$ and brine, with different saturations, and constant viscosity. They 
observe peak ratios of leakage rate to injection rate in the range of 0.2\% to 0.42\%. For the simulation presented in this work, 
the ratio is between 1\% and 10\%. Clearly, the assumptions of fully saturated, single phase flow do not accurately capture the 
leakage rate, but the results are meaningful in that they show, in general, an over-prediction of the leakage rate if a pressure 
dependent viscosity is not accounted for. As $\bar{\beta}$ increases, the amount $\mathrm{CO}_2$ that leaks into aquifer B is 
decreased significantly. Such results suggest that a pressure dependent viscosity model has a substantial effect on the predicted 
leakage rate.

\begin{table}[ht!]
\caption{$\mathrm{CO}_2$ leakage through an abandoned well: problem parameters}
\centering
\begin{tabular}{ll}
\hline
Parameter & Value  \\
\hline
$\bar{\alpha}_0$ aquifer A  & 1.0  \\
$\bar{\alpha}_0$ aquifer B  & 1.0  \\
$\bar{\alpha}_0$ inside the wells  & 100  \\
$\bar{p}_0$ aquifer A  & 1.0  \\
$\bar{p}_0$ aquifer B  & 0.95  \\
$\beta$ & $1.0 \times 10^{-7}$ \\
$\boldsymbol{v} \cdot \boldsymbol{n}$ at inflow & 0.262\\
\hline
\end{tabular}
\label{table:Probelm2Params}
\end{table}

%
\section{CONCLUDING REMARKS}
\label{Sec:S5_StabMDarcy_Conclusions}
We have considered a modification to Darcy equation by taking 
into account the dependence of viscosity on the pressure, which 
has been observed in many experiments on organic liquids. We 
have developed a new mixed stabilized formulation for the 
modified Darcy equation. We have also presented a numerical 
solution procedure to solve the resulting nonlinear equations. 
It has been shown numerically that equal-order interpolation for 
the velocity and pressure is stable under the proposed 
stabilized mixed formulation, which is not the case with 
the classical mixed formulation. 
Using representative problems, it has been shown that 
\begin{enumerate}[(a)]
\item the proposed stabilized mixed formulation performs 
  well, and 
\item the results predicted by the standard Darcy model 
  are qualitatively and quantitatively different from that 
  of the predictions based on the modified Darcy model. 
\end{enumerate}
It has been observed that the dependence of viscosity 
on the pressure drastically alters the pressure profile 
in the domain, and creates steep gradients near the 
boundary of the domain. Using a representative problem 
with relevance to geological carbon-dioxide sequestration, 
it has been shown that the standard Darcy model over-predicts 
the velocity of the fluid in the abandon well in comparison 
with the modified Darcy model. This prediction will serious 
consequences in designing the ceiling of the cap rock, 
which is one of the main mechanisms for the safety of 
a geological carbon-dioxide geosystem. Another important 
point to be noted is that the modified Darcy model predicts 
higher pressures than that of the Darcy model, which will 
have implications in modeling damage and fracture of the 
porous solid. 
 
As a part of future work on the numerical front, one can 
extend the current work to solve problems with much higher 
$\beta$ values by employing continuation-type methods (that 
is, to solve the problem for a high $\beta$ using information 
from the solution at a lower $\beta$). Another interesting 
future numerical work can be designing preconditioners for 
these kinds of nonlinear problems. A possible future work 
on the modeling front is to model the damage and fracture 
of the porous solid along with the flow aspects by taking 
into account the dependence of the viscosity on the pressure.

\section*{ACKNOWLEDGMENTS}
The first author (K.~B.~Nakshatrala) was supported in part by the 
Department of Energy through the Subsurface Biogeochemical 
Research (SBR) Program. The authors also acknowledge the 
financial support from Sandia National Laboratories through 
the Laboratory Directed Research and Development program. 
The opinions expressed in this paper are those of the authors 
and do not necessarily reflect that of the sponsors. The authors 
thank Pat Notz and Mario Martinez for their insightful discussions 
on this topic.

\bibliographystyle{unsrt}
\bibliography{../Master_References/Master_References,../Master_References/Books}
\clearpage
\newpage

\clearpage
\begin{figure}[htb!]
\centering
\includegraphics[scale=0.4]{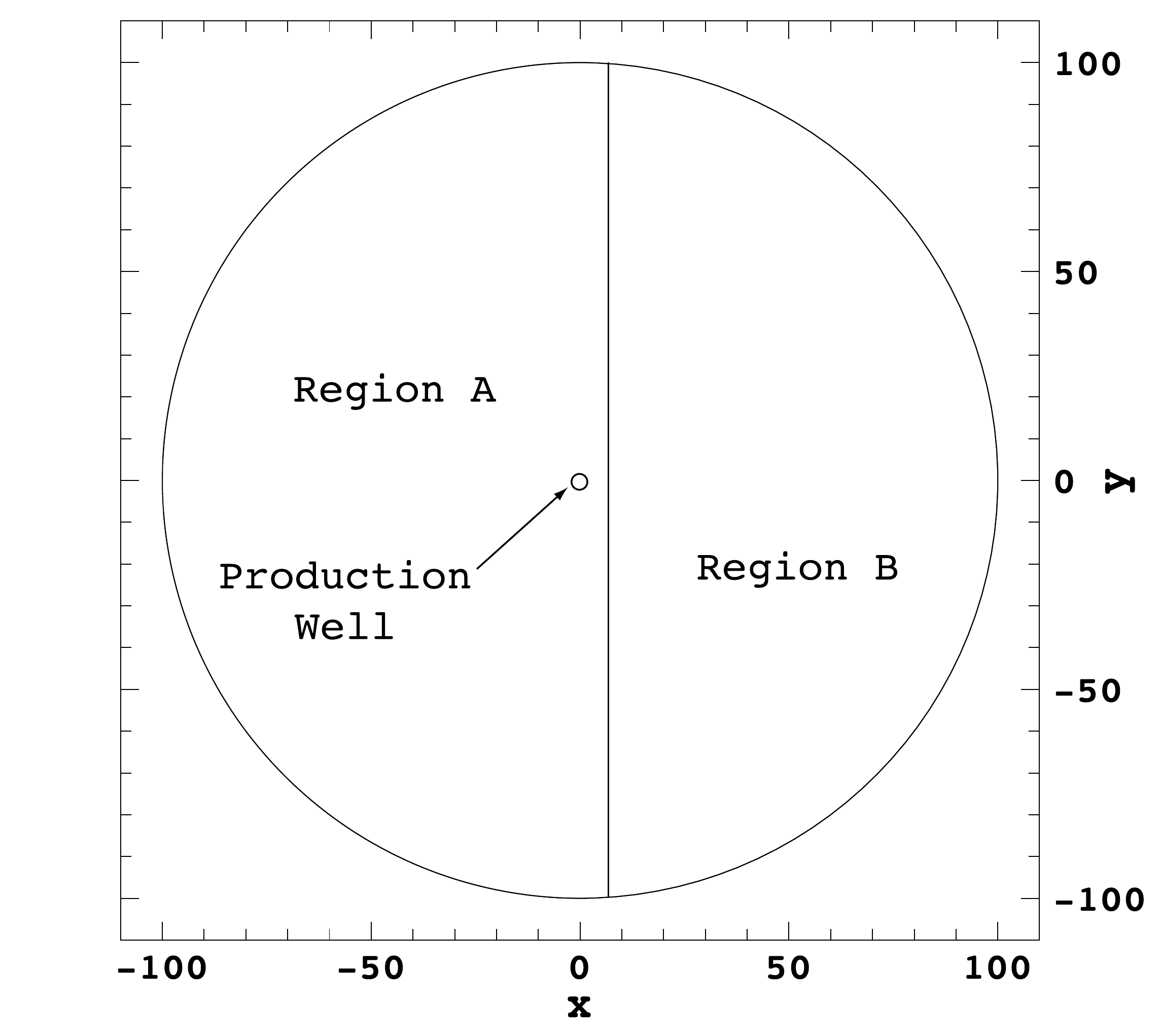}
\caption{Regions with different permeability: This figure shows the computational 
domain and the location of production well.} \label{fig:SunsetGeo}
\end{figure}

\begin{figure}[htb!]
\centering
\includegraphics[scale=0.4]{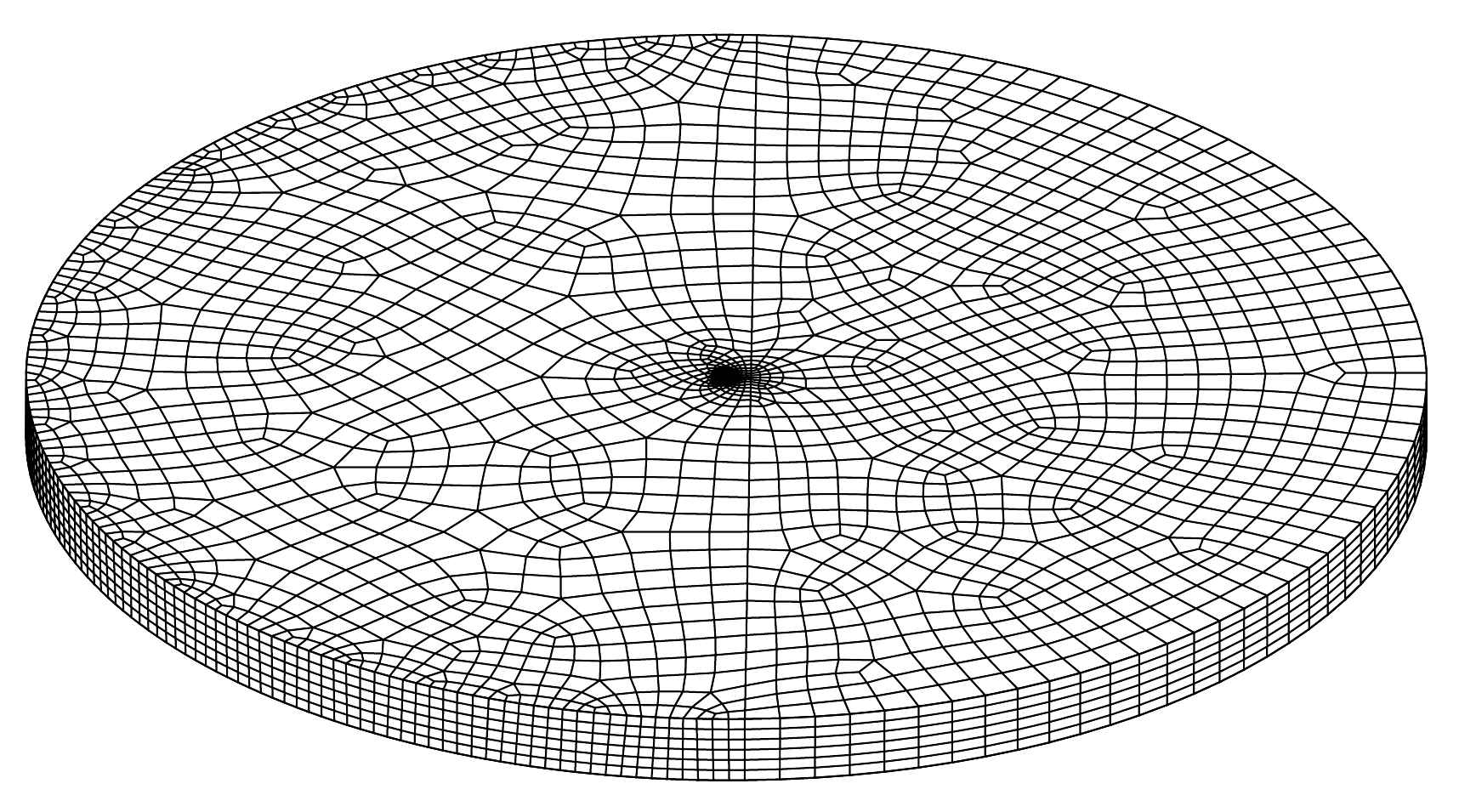}
\caption{Regions with different permeability: Three-dimensional finite element mesh 
using eight-node (linear) brick elements.} \label{fig:SunsetMesh}
\end{figure}

\end{document}